\documentclass[graybox]{svmult}

\usepackage{mathptmx}       % selects Times Roman as basic font
\usepackage{helvet}         % selects Helvetica as sans-serif font
\usepackage{courier}        % selects Courier as typewriter font
\usepackage{type1cm}        % activate if the above 3 fonts are
\usepackage[algo2e]{algorithm2e}
\usepackage{amsmath}
\usepackage{amssymb}
\usepackage{soul}

\usepackage{makeidx}         % allows index generation
\usepackage{graphicx}        % standard LaTeX graphics tool
\usepackage[caption=false]{subfig} % when including figure files

\usepackage{multicol}        % used for the two-column index
\usepackage[bottom]{footmisc}% places footnotes at page bottom
\usepackage{url}   
\usepackage{verbatim}   
\usepackage[utf8]{inputenc}
\usepackage{enumerate}
\usepackage{float}
\usepackage{todonotes}

\usepackage{tikz, pgfplots}
\usetikzlibrary{arrows,intersections,calc,patterns,shapes,arrows}
\makeatletter %new code
\pgfdeclarepatternformonly[\LineSpace,\tikz@pattern@color]{my north east lines}{\pgfqpoint{-1pt}{-1pt}}{\pgfqpoint{\LineSpace}{\LineSpace}}{\pgfqpoint{\LineSpace}{\LineSpace}}%
{
	\pgfsetcolor{\tikz@pattern@color} %new code
	\pgfsetlinewidth{0.4pt}
	\pgfpathmoveto{\pgfqpoint{0pt}{0pt}}
	\pgfpathlineto{\pgfqpoint{\LineSpace + 0.1pt}{\LineSpace + 0.1pt}}
	\pgfusepath{stroke}
}
\makeatother %new code
\newdimen\LineSpace
\tikzset{
	line space/.code={\LineSpace=#1},
	line space=3pt
}
\usepgfplotslibrary{fillbetween}
\usetikzlibrary{graphs} 
\usetikzlibrary{arrows}
\usetikzlibrary{matrix}
\pgfplotsset{compat=1.14}
\definecolor{darkblue}{RGB}{0,61,119}
\definecolor{lightblue}{RGB}{67,137,188}

\tikzstyle{block} = [rectangle, rounded corners, fill=darkblue!50, text width=7em, text centered, minimum height=5em, node distance=3cm]
\tikzstyle{cloud} = [ellipse,fill=lightblue!30, node distance=4cm, text width=5.5em, text centered, minimum height=5em]
\tikzstyle{ball} = [circle, fill=lightblue!30, text width=6em, text centered, minimum height=12em]

\newcommand{\Oad}{\mathcal{O}^{\text{ad}}} % set of admissible shapes
\newcommand{\R}{\mathbb{R}}	% real numbers

\newcommand{\pp}[2]{\frac{\partial #1}{\partial #2}} % partial#1/partial#2

\makeindex             % used for the subject index
\begin{document}
%%%%%%%%%%%%%%%%%%%%%% begin of template for preprint %%%%

%{

%\setcounter{page}{0}

%\thispagestyle{empty}

%\setlength{\parindent}{0cm}

%\setlength{\parskip}{0cm}

%\begin{center}

%\mbox{}

%\vspace{-2cm}

%\large

%\includegraphics[width=5cm]{logo_imacm.pdf}

%\vspace{2cm}

%Bergische Universit\"at Wuppertal

%\vspace{0.5cm}

%Fachbereich Mathematik und Naturwissenschaften

%\vspace{0.5cm}

%Institute of Mathematical Modelling, Analysis and Computational Mathematics (IMACM)

%\vspace{1cm}

%Preprint BUW-IMACM
%% Preprint-Identification
%20/02

%%{\small Updated version of Preprint No. 11/01}

%\vspace{1.5cm}

%% Name(s) of author(s)
%Jan Backhaus, Matthias Bolten, Onur Tanil Doganay, Matthias Ehrhardt, Benedikt Engel, Christian Frey, Hanno Gottschalk, Michael G\"unther, Camilla Hahn, Jens J\"aschke, Peter Jaksch, Kathrin Klamroth, Alexander Liefke, Daniel Luft, Lucas M\"ade, Vincent Marciniak, Marco Reese, Johanna Schultes, Volker Schulz, Sebastian Schmitz, Johannes Steiner and Michael Stiglmayr

%\vspace{1cm}

%{\bf \Large
%% Title of preprint

%GivEn -- Shape Optimization for Gas Turbines in Volatile Energy Networks

%}

%\vspace{1cm}

%% Date of publication (only month and year)
%February 2020

%\vspace{1cm}

%{\sf
%http://www.math.uni-wuppertal.de 
%}
%\end{center}

%}

%\clearpage
%%%%%%%%%%%%%%%%%%%%%% end of template for preprint %%%%

\title*{GivEn -- Shape Optimization for Gas Turbines in Volatile Energy Networks}
\author{
Jan Backhaus
\and Matthias Bolten
\and Onur Tanil Doganay
\and Matthias Ehrhardt
\and Benedikt Engel
\and Christian Frey
\and Hanno Gottschalk
\and Michael G\"unther
\and Camilla Hahn
\and Jens J\"aschke
\and Peter Jaksch
\and Kathrin Klamroth
\and Alexander Liefke
\and Daniel Luft
\and Lucas M\"ade
\and Vincent Marciniak
\and Marco Reese
\and Johanna Schultes
\and Volker Schulz
\and Sebastian Schmitz 
\and Johannes Steiner
\and Michael Stiglmayr}
\authorrunning{The GivEn Consortium} 

\institute{
Jan Backhaus, Christian Frey \at Institute of Propulsion Technology,
German Aerospace Center (DLR), 51147 K\"oln, Germany,
\email{{jan.backhaus,christian.frey}@dlr.de}
\and Matthias Bolten, Onur Tanil Doganay, 
Matthias Ehrhardt, Hanno Gottschalk, 
Michael G\"unther, Camilla Hahn, Jens J\"aschke, Kathrin Klamroth, Marco Reese, 
Johanna Schultes and Michael Stiglmayr
\at Bergische Universität Wuppertal,
Fakultät f\"ur Mathematik und Naturwissenschaften, IMACM, Gaußstrasse 20, 42119 Wuppertal, Germany,
\email{$\{$bolten,doganay,ehrhardt,guenther,hgotsch,chahn,jaeschke, klamroth,reese,jschultes,stiglmayr$\}$@uni-wuppertal.de}
\and Vincent Marciniak, Alexander Liefke and Peter Jaksch \at
Siemens AG, Power and Gas, Common Technical Tools
Mellinghoffer Str. 55, 45473 M\"ulheim an der Ruhr, Germany,
\email{$\{$vincent.marciniak,alexander.liefke,peter.jaksch$\}$@siemens.com}
\and Daniel Luft and Volker Schulz \at
Universität Trier, Fachbereich IV,
Research Group on PDE-Constrained Optimization,
54296 Trier, Germany,
\email{$\{$luft,volker.schulz$\}$@uni-trier.de}
\and Lucas M\"ade, Johannes Steiner and Sebastian Schmitz \at 
Siemens Gas and Power GmbH \& Co. KG, Probabilistic Design,
GP PGO TI TEC PRD,
Huttenstr. 12, 10553 Berlin, Germany,
\email{$\{$lucas.maede,johannes.steiner,schmitz.sebastian$\}$@siemens.com}
\and Benedikt Engel \at 
University of Nottingham, Gasturbine and Transmission Research Center (G2TRC),
NG72RD Nottingham, United Kingdom,
\email{engel.benedikt@nottingham.ac.uk}
}
%
% Use the package "url.sty" to avoid problems with special characters
\maketitle

%%%%%%%%%%%%%%%%%%%%%%%%%%%%%%%%%%%%%%% ABSTRACT %%%%%%%%%%%%
\abstract*{
This paper describes the project GivEn that develops a novel multicriteria optimization process for gas turbine blades and vanes using modern ''adjoint'' shape optimization algorithms. 
Given the many start and shut-down processes of gas power plants in volatile energy grids, besides optimizing gas turbine geometries for efficiency, the durability understood as minimization of the probability of failure is a design objective of increasing importance. 
We also describe the underlying coupling structure of the multiphysical simulations and use modern, gradient based multicriteria optimization procedures to enhance the exploration of Pareto-optimal solutions.
}

\abstract{
This paper describes the project GivEn that develops a novel multicriteria optimization process for gas turbine blades and vanes using modern ''adjoint'' shape optimization algorithms. 
Given the many start and shut-down processes of gas power plants in volatile energy grids, besides optimizing gas turbine geometries for efficiency, the durability understood as minimization of the probability of failure is a design objective of increasing importance. 
We also describe the underlying coupling structure of the multiphysical simulations and use modern, gradient based multicriteria optimization procedures to enhance the exploration of Pareto-optimal solutions.
}

%\abstract{
%The project GivEn develops a novel multicriteria optimization process for gas turbine blades and vanes using modern ``adjoint'' algorithms from the mathematical discipline of shape optimization.  Due to frequent start and shut-down processes of gas power plants in volatile energy grids,  the durability understood as minimization of the probability of failure is a design objective of increasing importance besides optimizing gas turbine geometries for aerodynamically efficiency. We describe the underlying coupling structure of the multiphysical simulation and use modern, gradient based multiobjective optimization to enhance the exploration of Pareto-optimal solutions.
%}

%%%%%%%%%%%%%%%%%%%%%%%%%%%%%%%%%%%%%%%%%%%%%%% 5 pages ( M. Ehrhardt)
\section{Introduction}\label{sec:intro}
\begin{figure}[htb!]
\centering
\resizebox{0.85\textwidth}{!}{\begin{tikzpicture}[scale=1, auto, align = center,font=\sffamily]
% circle radius
\def \klR {4.5cm}
\def \grR {8cm}
\def \a {11} % arc shorten for (eff) and (temp_distr)
\def \b {18} % arc shorten for (BC) and (PoF)
\def \c {18} % arc shorten for one part of the blocks (more in the middle)
\def \d {14} % arc shorten for one part of the blocks (more away from the middle)
\def \usertextsize {\small}

% blocks and cloud elements
\node[block](HC) at ({360/4 * 0.5}:\klR) {\usertextsize{heat conduction}};
\node[block](CC) at ({360/4 * 0.5}:\grR) {\usertextsize cooling channels};
\node[block](CFD) at ({360/4 * 1.5}:\klR) {\usertextsize CFD-model};
\node[block](MOP) at ({360/4 * 2.5}:\klR) {\usertextsize gradient-based multiobjective optimization};
\node[block](TMG) at ({360/4 * 3.5}:\klR) {\usertextsize thermo-mechanical equation};

\node[cloud](temp_distr) at ({360/4 * 0}:\klR) {\usertextsize temperature distribution};
\node[cloud](BC) at ({360/4 * 1}:\klR) {\usertextsize boundary  conditions};
\node[cloud](eff) at ({360/4 * 2}:\klR) {\usertextsize aerodynamic objective};
\node[cloud](LCF) at ({360/4 * 3}:\klR) {\usertextsize LCF objective};

\node[ball](update) at (0:0) {\usertextsize \textbf{geometry update of the component}};

% arcs
\draw[->, >=latex, line width=1pt] ({360/4 * 4.0-\a}:\klR) arc ({360/4 * 4.0-\a}:{360/4 * 3.5+\d}:\klR); % temp distr -> TM eq. 
\draw[->, >=latex, line width=1pt] ({360/4 * 3.5-\c}:\klR) arc ({360/4 * 3.5-\c}:{360/4 * 3.0+\b}:\klR); % TM eq. -> PoF
\draw[->, >=latex, line width=1pt] ({360/4 * 3.0-\b}:\klR) arc ({360/4 * 3.0-\b}:{360/4 * 2.5+\c}:\klR); % PoF -> MOP
\draw[->, >=latex, line width=1pt] ({360/4 * 2.0+\a}:\klR) arc ({360/4 * 2.0+\a}:{360/4 * 2.5-\d}:\klR); % MOP <- eff 
\draw[->, >=latex, line width=1pt] ({360/4 * 1.5+\d}:\klR) arc ({360/4 * 1.5+\d}:{360/4 * 2.0-\a}:\klR); % eff <- CFD 
\draw[->, >=latex, line width=1pt] ({360/4 * 1.5-\c}:\klR) arc ({360/4 * 1.5-\c}:{360/4 * 1.0+\b}:\klR); % CFD -> BC
\draw[->, >=latex, line width=1pt] ({360/4 * 1.0-\b}:\klR) arc ({360/4 * 1.0-\b}:{360/4 * 0.5+\c}:\klR); % BC -> HC
\draw[->, >=latex, line width=1pt] ({360/4 * 0.5-\d}:\klR) arc ({360/4 * 0.5-\d}:{360/4 * 0.0+\a}:\klR); % HC -> temp distr.  

\draw[->, >=latex, line width=1pt, shorten >=2pt,shorten <=2pt] (HC) to [bend right=15] (CC);
\draw[->, >=latex, line width=1pt, shorten >=2pt,shorten <=2pt] (CC) to [bend right=15] (HC);

\draw[->, >=latex, line width=1pt, shorten >=2pt,shorten <=2pt] (update) to (CFD.south east);
\draw[->, >=latex, line width=1pt, shorten >=2pt,shorten <=2pt] (update) to (HC.south west);
\draw[->, >=latex, line width=1pt, shorten >=2pt,shorten <=2pt] (update) to (TMG.north west);
\draw[->, >=latex, line width=1pt, shorten >=2pt,shorten <=2pt] (MOP.north east) to (update);

% % check if the arcs are lying on the same circle 
% \draw (0,0) circle (\klR);
\end{tikzpicture}}
\caption{Information flow and dependencies between project parts.\label{fig:model}}
\end{figure}
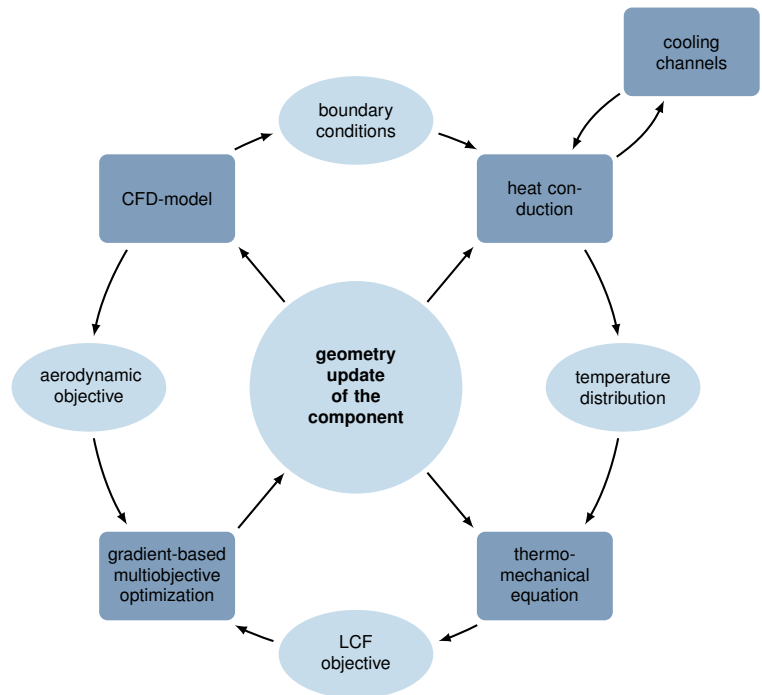

The diverse applications of gas turbines in the context of the energy system transformation, such as backup power plants or hydrogen turbines, go hand in hand with specific design requirements, in particular with regard to the efficiency of energy conversion and the reliability and flexibility of operation. 
These different requirements are intensively related to the coupled fluid dynamic simulation and the structural mechanical fatigue calculation.
The use of integrated, multi-physical tool chains and optimization software therefore plays an important role in gas turbine design. 
This joint project links six different simulations -- fluid dynamics, laminar convective heat transfer, 1D flux networks and turbulent convective heat transfer, heat conduction, structural mechanics, probabilistic modelling of material fatigue -- which are computed on a complex turbo geometry.
These simulations are coupled in the multi-objective shape optimization process. See Fig.~\ref{fig:model} for a schematic illustration of the multi-physical simulation/optimization cycle.

The challenge for the GivEn project thereby is to adjoin a highly multi-physical simulation chain with continuous coupling, to determine form gradients and form Hessians with respect to different objectives, and to make it usable in the multicriteria optimization for the turbine design process. 

Coupled multi-physics simulation is an ongoing topic in turbo-machinery. 
For recent surveys of these fluid dynamics and heat transfer topics, see e.g.\ \cite{sultanian2018gas,tyacke2019turbomachinery}.
The important topic of turbo-machinery life calculation is often treated separately and from a materials science point of view, see e.g.\ \cite{baeker08,cui2002state,rao2008lifing}. 
In contrast to the traditional separation of the mechanical and the fluid dynamics properties, 
the approach we follow in GivEn preserves a holistic viewpoint.

The algorithmic optimization of turbo-machinery components by now has a long history. 
While in the beginning genetic algorithms were used predominantly, in recent times data driven methods like Gaussian processes or (deep) neural networks predominate \cite{adams2009dakota,dorfner2007axis,siller2009automated}. 
The strength of such procedures lie in global search approaches.   
As an alternative, gradient based optimization using adjoint equations are seen as a highly effective local search method \cite{frey2009,giles2000introduction}, see also \cite{li2017} for a recent review including a comparison of the methods and \cite{backhaus2017application} for bringing the data driven and the adjoint world together using gradient enhanced Gaussian processes \cite{adams2009dakota,forrester2008engineering}.

When combining the challenge of multi-physics and multi-criteria optimization, it would be desirable to treat mechanical and fluid dynamic aspects of turbo-machinery design on the same footing.   
A necessary prerequisite for this is the probabilistic modelling of the mechanisms of material damage, 
as this enables the application of the adjoint method \cite{bolten2015,gottschalk2015,gottsch,gottschalk2018adjoint,gottschalk2014,gottschalk2014b,SchmitzSeibel_2013,MaedeSchmitz_2018}. 
This is not possible with a deterministic calculation of the lifetime of the weakest point, as taking the minimum over all points on the component is a non differential operation.

The GivEn consortium exploits these new opportunities for multi-criteria and multiphysics optimization. 
It brings together a leading original equipment manufacturer (Siemens Power and Gas), technology developing institutions 
(German Aero Space Center (DLR) and Siemens CT) as well as researchers from academia (Universities of Trier and Wuppertal). 
Since 2017 this consortium addresses the challenges described in a joint research effort funded by the BMBF under the funding scheme "mathematics for innovation". 
With the present article, we review the research done so far and give an outlook on future research efforts.

%%%%%%%%%%%%%%%%%%%%%%%%%%%%%%%%%%% structure of paper ....
This paper is organized as follows. 
In Section~\ref{sec:projects} we describe our research work on the different physical domains including the usage of adjoint equations, improved shape gradients and gradient based multi-criteria optimization. 
Following the design scheme outlined in Figure~\ref{fig:model}, we start with aerodynamic shape optimization in Section~\ref{sec:21} using modern mesh morphing based on the Steklov-Poincar\'e definition of shape gradients \cite{schulz2014,Schulz2016b,Schulz2016}, then proceed to heat transfer and the thermal loop in Section~\ref{sec:22}. 
Section~\ref{sec:23} includes related probabilistic failure mechanisms. 
The model range from empirical models based on Weibull analysis and point processes to elaborate multi scale models. 
Section~\ref{ShapeOpStruc} presents shape optimization methods that are based on the probability of failure and develops a highly efficient computational framework based on conformal finite elements.  
Section~\ref{sec:25} presents novel fundamental results on the existence of Pareto fronts in shape optimization along with algorithmic developments in multi-criteria gradient based shape optimization including scalarization, bi-criteria gradient descent and gradient enhanced Gaussian processes. 
% Next, in Section~\ref{sec:26} a probabilistic model for LCF is investigated.

In Section~\ref{sec:applications} we describe the industrial perspective from the standpoint of the DLR and Siemens energy. 
While in Section~\ref{sec:31} the DLR gives a description of the interfaces with and the possible impact to the DLR's own R\&D roadmap, Siemens Power \& Gas in Section~\ref{sec:32} relates adjoint based multi-criteria optimization with adjoint multi-criteria tolerance design and presents an application on real world geo\-metries of 102 casted and scanned turbine vane geometries. 

Let us note that this work is based on the papers \cite{Bolten2019,Doganay2019,engel2019probabilistic,gottsch,gottschalk2018adjoint,liefke2020,liefke2019,liefke2018,luft2019efficient,MaedeSchmitz_2018} that have been published with (partial) funding by the GivEn consortium so far. 
As this report is written after about half of the funding period of the project, we also give comments on future research plans within GivEn and beyond.

%%%%%%%%%%%%%%%%%%%%%%%%%%%%%%%%%%%%%%%%%%%%%
\section{Areas of Mathematical Research and Algorithmic Development}\label{sec:projects}
The project GivEn researches the multiobjective free-form optimization of turbo geometries. 
For this purpose, the thermal and mechanical stress of the turbine blades and their aerodynamic behavior must be modelled, simulated and optimised.
In the following we describe the components of the multiphysical simulation and optimization, namely aerodynamic shape optimization, heat transfer and thermal loop, probabilistic objective functionals for cyclic fatigue, shape optimization for probabilistic structure mechanics, multiobjective optimization, and probabilistic material science.

%%%%%%%%%%%%%%%%%%%%%%%%%%%%%%%%%%%%%%%%%% WP1 (PI Schulz, U Trier)
\subsection{Aerodynamic Shape Optimization}\label{sec:21}
Shape optimization is an active research field in mathematics. 
Very general basic work on shape calculus can be found in \cite{Sokolowski92, Haslinger03, Eppler07}. 
Aerodynamic investigations can be found in \cite{SchmidtThesis10, Schmidt10}.  
New approaches understand shape optimization as the optimization on shape manifolds \cite{schulz2014, WelkerThesis16} and thus enable a theoretical framework that can be put to good practical use, while at the same time leading to mathematical challenges, as no natural vector space structure is given. Otherwise, applications usually use finite dimensional parameterizations of the form, which severely limits the space of allowed shapes. 

In the shape space setting, the use of volume formulations has been shown in combination with form metrics of the Steklov-Poincar\'e type \cite{Schulz2016, Schulz2016b} were shown to be  numerically very advantageous, since the volume formulation in comparison to the formally equivalent boundary formulation for canonical discretizations have better approximation properties and also weaker smoothness requirements of the functions involved. 
Additionally the Steklov-Poincar\'e type metrics require a free combination of volume and boundary formulations together with an inherently good approximation of the Shape-Hessian operators.

% trasor package introduction
In order to exploit these theoretical advances for industrial applicability meeting high-end standards, 
the TRASOR (TRACE Shape Optimization Routine) software package for non-parametric shape optimization routines has been created.
This software package is built on several solver bundles connected by an interface in Python~2.7 and 3.5. 
One major package bundle provided by the DLR and incorporated in TRASOR is TRACE~9.2, which is an interior flow simulator (cf.\ \cite{becker2010}).
%\textcolor{red}{(@Jan: gibts da ein standard Paper/Link den wir dazu zitieren können?)}
%
%---------------------------
%The first major package bundle developed and provided by the DLR for interior flow simulation in turbomachinery, and incorporated in TRASOR consists of \textcolor{red}{Versionen von packages; wieviel Detail?} TRACE 9.2.426 [49, 14], POST 9.2.426 and PREP 9.2.426 .
%---------------------------
%
The TRASOR software incorporates shape gradient representations using Steklov-Poincar\'{e}-metrics (cf.\ \cite{Schulz2016, Schulz2016b}) based on shape sensitivities derived by automatic differentiation provided by adjointTRACE~\cite{Sagebaum2017, backhaus2017application}.

%% Fenics intro kurz + Stek Poincare
TRASOR also interfaces with FEniCS 2017.2.0 \cite{AlnaesBlechta2015a, LoggMardalEtAl2012a}, which is a Python based finite element software utilizing several sub-modules, such as the Unified Form Language (UFL \cite{Alnaes2012a}), Automated Finite Element Computing (DOLFIN \cite{LoggWells2010a, LoggWellsEtAl2012a}) and PETSc \cite{balay2019petsc} as a linear algebra backend, in order to solve  differential equations based weak formulations. 
Various solver options, including CG, GMRES, PETCs's built in LU solver and preconditioning using incomplete
LU and Cholesky, SOR or algebraic multigrid methods are available in FEniCS and thus applicable in TRASOR. 
FEniCS/PETSc also offers the possibility to parallelize finite element solving, making the Steklov-Poincar\'e gradient calculation scalable in processor number. 

\medskip
\noindent Features of the software package TRASOR include
%\todo{H: Do not use bullet points, as this is the only place in the article where bullets are used.} done ME:
\renewcommand*\labelitemi{\normalfont\textendash}
\begin{itemize}
	\item automatic file generation and management for TRACE and adjointTRACE
	\item interface between TRACE and FEniCS, including automatic FEniCS mesh gene\-ration from .cgns files
	\item steepest descent optimization using TRACE intern gradients
	\item steepest descent optimization using Steklov-Poincar\'{e} gradients calculated in FEniCS
	\item target parameter selection for various parameters found in TRACE, including all parameters listed in \cite{TRACEV}
	\item generation of .pvd and .vtu files of gradients, sensitivities, meshes and flow simu\-lation data for visual post processing
\end{itemize}

TRASOR features are tested on the low-pressure turbine cascade T106A designed by MTU Aero Engines (cf.\ \cite{hoheisel1986influence}).
The algorithm using  Steklov-Poincar\'e gradients is outlined in Algorithm~\ref{algo:TRASOR}.

\begin{figure}
\centering
\begin{tikzpicture}[thick,scale=0.7, every node/.style={scale=0.7}]
	\node[block] (a){TRACE.cgns};
	\node[block, right = 0.8cm of a] (b){TRACE.dat};
	\node[block, right = 0.8cm of b] (c){FEniCS.msh};
	\node[block, right = 0.8cm of c] (d){FEniCS.xdmf/-.h5};
	\node[block, right = 0.8cm of d] (e){FEniCS mesh};
	
	\draw[thick,->] (a.east) -- (b.west) node[midway,above] {POST};
	\draw[thick,->] (b.east) -- (c.west) node[midway,above] {TRASOR};
	\draw[thick,->] (c.east) -- (d.west) node[midway,above] {meshio};
	\draw[thick,->] (d.east) -- (e.west) node[midway,above] {DOLFIN};
\end{tikzpicture}
	\caption{TRACE to FEniCS Pipeline}
	\label{fig::TRACE2FEniCS}
\end{figure}
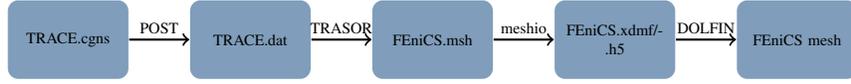

In order to exploit FEniCS it is necessary to create an unstructured computational mesh with vertices prescribed by TRACE. 
As FEniCS 2017.2.0 is not fully capable of supporting hexahedral and quadrilateral elements (this should be available with FEniCS 2020), hexahedral and quadrilateral elements used in TRACE are partitioned to conforming tetrahedral and triangular elements respecting the structured TRACE mesh. The conversion process including the data formats for TRACE to FEniCS mesh conversion are depicted in Fig.~\ref{fig::TRACE2FEniCS} (cf.\ \cite{meshio2019, LoggWellsEtAl2012a})
%\begin{align}\label{eq:TRASOR_Meshformats}
%\begin{split}
	%&\text{TRACE.cgns}
	 %\overset{\text{POST}}{\longrightarrow} \text{TRACE.dat}  %\overset{\text{TRASOR}}{\longrightarrow} \\ &\text{FEniCS.msh}  %\overset{\text{meshio}}{\longrightarrow} \text{FEniCS.xdmf/-.h5} %\overset{\text{DOLFIN}}{\longrightarrow} \text{FEniCS mesh}
%\end{split}
%\end{align}

For representing the TRACE generated mesh sensitivities $D_{ad}J(\Omega_{\text{ext},k})$ as a Steklov-Poincar\'e gradient a sufficient metric has to be chosen. 
According to \cite{Schulz2016}, we implemented the following linear elasticity model 

\begin{align}\label{eq:LinearElasticityTRACE}
\begin{split}
\int_{\Omega_{\text{ext},k}}\sigma\bigl(\nabla^{StP}J(\Omega_{\text{ext},k})\bigr):\epsilon(V)\,dx &= D_{ad}J(\Omega_{\text{ext},k})[V] \quad \forall\,V\in H^1_0(\Omega_{\text{ext},k}, \R^d) \\
\nabla^{StP}J(\Omega_{\text{ext},k}) &= 0 \qquad \qquad \qquad \text{ on } \Gamma_{\rm Inlet/Outlet} \\
\sigma(V) &= \lambda \operatorname{Tr}\bigl(\epsilon(V)\bigr)I + 2\mu \epsilon(V) \\
\epsilon(V) &= \frac{1}{2}(\nabla V + \nabla V^\top),
\end{split}
\end{align}
where $\lambda\in\mathbb{R}$, $\mu\in\mathbb{R}_+$ are the so called Lam\'e parameters. 
If $D$ is the entire duct including the shape of the turbine blade $\Omega_{k}$ at iteration $k$ of the shape optimization procedure, $\Omega_{\text{ext},k}=D\setminus \Omega_{k}$ is the external computational domain where the fluid dynamics takes place. 
Continuous Galerkin type elements of order one are used for target and test spaces in the FEniCS subroutine conducting the shape gradient calculation. 

\begin{figure}[t]
	\centering
	\subfloat[Steklov-Poincar\'e gradient]{\includegraphics[width=0.485\linewidth]{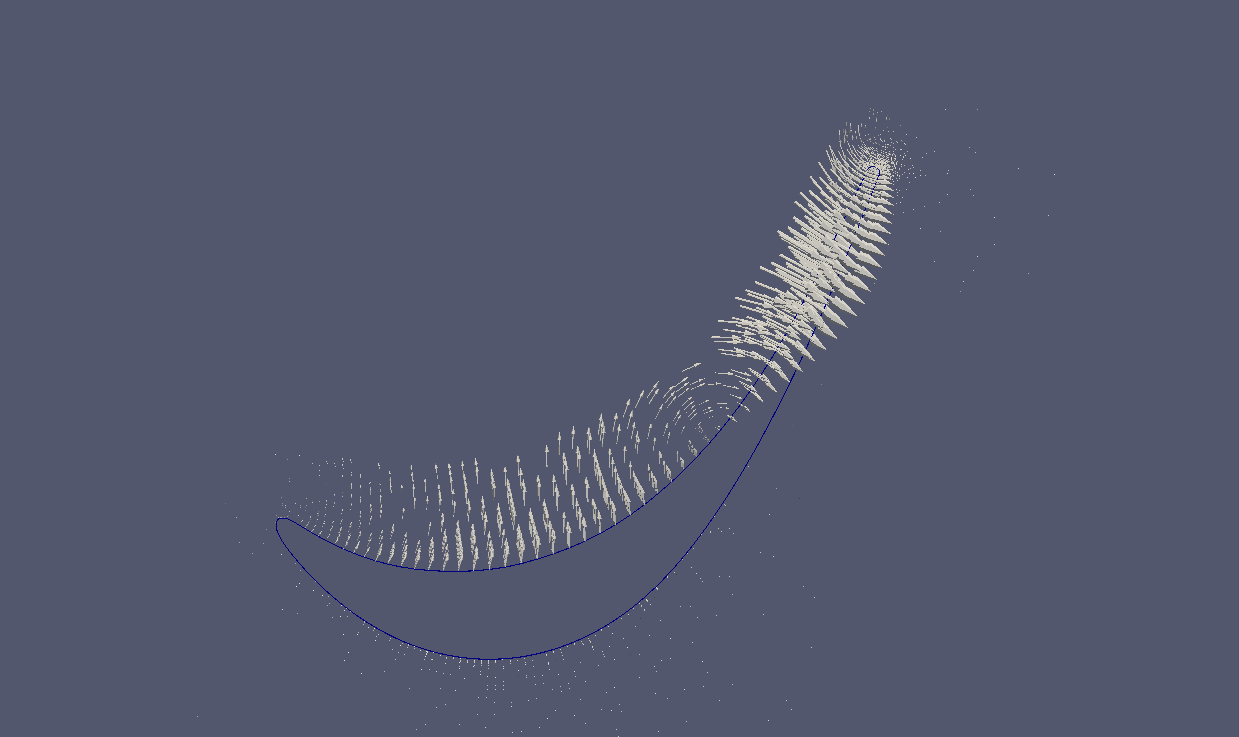}}
	\quad
	\subfloat[TRACE Mesh Smoothing gradient]{\includegraphics[width=0.485\linewidth]{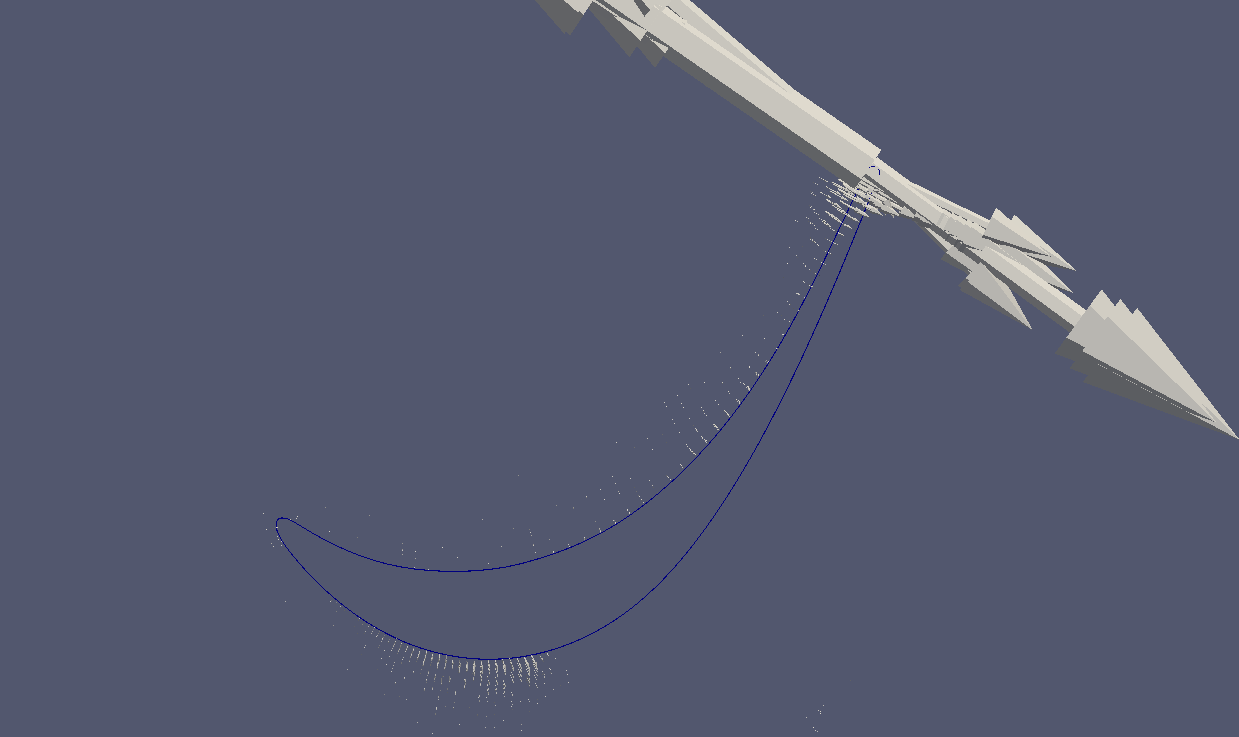}}
\\
	\subfloat[FEniCS computational mesh of the T106A]{\includegraphics[width=0.485\linewidth]{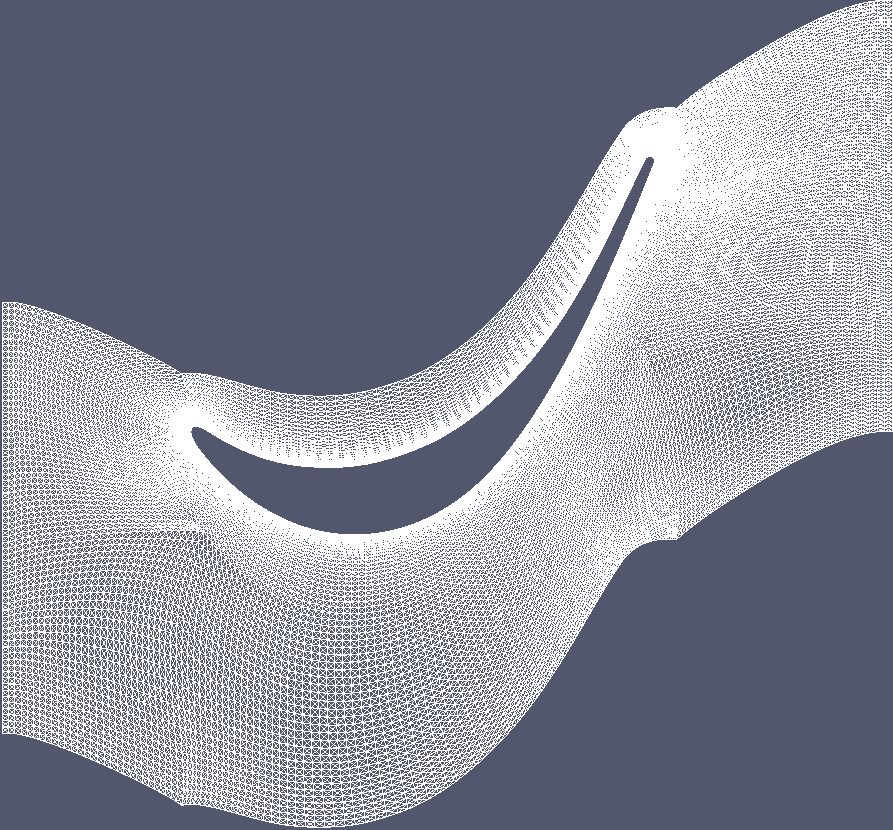}}
	\caption{Comparison of Steklov-Poincar\'e gradient (upper left) with TRACE Mesh Smoothing gradient (upper right) on a FEniCS computational mesh of the T106A (lower center)}
	\label{fig:TRACE}
	\vspace*{1cm}
\end{figure}

An exemplary comparison of a Steklov-Poincar\'e gradient calculated by solving the linear elasticity system~\eqref{eq:LinearElasticityTRACE} with Lam\'e parameters $\lambda\equiv 0$ and constant $\mu >0$, and a TRACE gradient, which is generated by solving a linear elasticity mesh smoothing system with Dirichlet boundaries being the lattice sensitivities $D_{ad}J(\Omega_{\text{ext},k})$,
for the isentropic total pressure loss coefficient in relative frame of reference based on dynamic pressure is portrayed in Fig.~\ref{fig:TRACE}.
We can see additional gain of regularity in the gradient through Steklov-Poincar\'e representation, in particular the pronounced rise in sensitivity at the trailing edge is handled by redistributing sensitivities at the pressure side in a smooth manner, thus guaranteeing better stability of the mesh morphing routine.

An in-depth comparison of shape optimization routines involving both types of gradient representation will be subject of a follow-up study. 
Further, a Steklov-Poincar\'e gradient representation using different bilinear forms matching the shape Hessian of the RANS flow and the target at hand are object of further studies, which might open new possibilities with superior convergence and mesh stability behavior. 

%% Bilder + erste Resultate

The following Algorithm~\ref{algo:TRASOR} is a prototype for a shape optimization problem, including the Steklov-Poincar\'e gradient representation. 
% Algorithm for TRASOR StekPoincare
\LinesNumbered
\DontPrintSemicolon
\begin{algorithm2e}\caption{TRASOR algorithm using Steklov-Poincar\'{e} gradients}
	\label{algo:TRASOR}
	\SetAlgoRefName{algo:TRASOR}
	Set flow parameters in TRACE.cgns, optimization parameters and targets in TRASOR.py  \;
	Build TRASOR file architecture \;
	Assemble and load FEniCS data from TRACE.cgns \;
	\While{$\Vert \nabla^{StP}J(\Omega_{\text{ext},k} \Vert > \varepsilon_{\text{\emph{shape}}}$}{
		Flow simulation and (AD) checkpoint creation using TRACE \;
		Calculate mesh sensitivities by automatic differentiation using adjointTRACE \;
		Pass mesh sensitivities to FEniCS setup \;
		Generate Steklov-Poincar\'{e} gradient in FEniCS: \;
		\Indp
			Calculate Lam\'{e}-Parameters \;
			Solve linear elasticity problem \eqref{eq:LinearElasticityTRACE} \;
		\Indm
		Extract target and flow values to update/ create protocols and .pvd/.vtu files \;
		Deform FEniCS mesh using FEniCS Steklov-Poincar\'e gradient and ALE (Arbitrary Lagrangian-Eulerian) \;
		Create TRACE\_deformation.dat files from FEniCS Steklov-Poincar\'e gradient \;
		Deform TRACE mesh using PREP	
	}
\end{algorithm2e} 

\clearpage % instead of 'newpage'

%%%%%%%%%%%%%%%%%%%%%%%%%% WP2 (PIs Ehrhardt and Günther)
% \newpage
\subsection{Heat Transfer and the Thermal Loop}\label{sec:22}

The numerical simulation of coupled differential equation systems is a challenging topic.
The difficulty lies in the fact that the (P)DEs involved may differ in type and also in order, and thus require different types and quantities of boundary conditions. \cite{bartel2018}
The question of the correct coupling is closely related to the construction of so-called transparent boundary conditions, which are based on the coupling of interior and exterior solutions.

The numerical simulation of coupled differential equation systems by means of co-simulation has the innate advantage that one can choose the optimal solver for each sub-system, for example by employing pre-existing simulation software. 
Most of the work done in this field concerns transient, i.e.\ time dependent, problems. 
In our case, however, we are interested in steady state systems which rarely get special attention in current research. 

Our model problem arises from the heat flow in a gas turbine blade.
Since higher combustion temperatures result in better efficiency \cite{saravanamuttoo2001gas}, engineering always strives for means to achieve these. 
However, this is limited by the material properties of the turbine blade, especially its melting point. 
One way to mitigate this, is by cooling the blade from the inside. 
This is done by blowing air through small cooling ducts. 
These ducts have a complex geometry to increase turbulence of the airflow and maximize heat transfer from the blade to the relatively cool air. For an overview, see for example \cite{han2004}.
Due to the small length-scales and high turbulence, regular fluid dynamics simulation techniques are infeasible for the simulation of the airflow within the ducts. 
Instead, they are modeled as a one-dimensional flow with parametric models for friction and heat transfer, similar to the work in \cite{meitner1990} and \cite{sultanian2018gas}.

\begin{equation}\label{eq:cooling_differential}
    \begin{split}
    w \pp{v}{x} &= A \pp{p}{x} + \frac{A}{2 D_h} f \rho v^2 + A \rho \omega^2 r \pp{r}{x} \\
    \pp{(v T)}{x} &= S \\
    \rho &= \frac{w}{v A}\\
    p &= \rho R_s T
\end{split}
\end{equation}

Here, $w$ is the mass flow through the channel which is assumed constant (i.e.\ only one inlet and outlet), $v$ is the fluid velocity, $\rho$ is the fluid's density and $p$ and $T$ denote the pressure and Temperature of the fluid, respectively. 
$A$ is the cross-sectional area and $D_h$ the hydraulic diameter of the channel, with $f$ as the Fanning friction factor. $\omega$ and $r$ are only relevant in the rotating case and denote the angular velocity and distance from the axis of rotation. 
$R_s$ is the specific gas constant of the fluid and $S$ denotes the heat source term from the heat flux through the channel walls.

The system~\eqref{eq:cooling_differential} takes the form of a DAE, but can be transformed into a system of ODEs by some simple variable substitutions.
Since the physical motivation behind the terms is easier to understand in the DAE form, this is omitted here. 

The heat conduction within the blade material is modeled by a PDE. 
In the transient case, this would be a heat equation.
In the stationary case, it is given by a Laplace equation. 
The heat transfer across the boundary is given by Robin boundary conditions, that prescribe a heat flux across the boundary depending on the temperature difference between ``inside" and ``outside". 
Temperature in \eqref{eq:conduction_differential} is denoted by $U$ to signify that it is mathematically a different entity than the temperature in the cooling channel, denoted by $T$ in \eqref{eq:cooling_differential}. $k$ is the thermal conductivity of the blade metal, while $h_{int}$ and $h_{ext}$ are the heat transfer coefficients of the internal and external boundary.

\begin{equation}\label{eq:conduction_differential}
    \begin{aligned}
        k \nabla^2 U &= 0 && \text{on~} \Omega\\
        -k \pp{U}{n} &= h_{\rm int} (U-U_{\rm int})  && \text{on~} \partial \Omega_{\rm int} \\
        -k \pp{U}{n} &= h_{\rm ext} (U-U_{\rm ext})  && \text{on~} \partial \Omega_{\rm ext}
    \end{aligned}
\end{equation}

The coupling between equations \eqref{eq:cooling_differential} and \eqref{eq:conduction_differential} is realized via the boundary condition, more specifically the internal boundary temperature $U_{int}$ as a function of $T$, on the conduction side and the source term $S$ in the cooling duct equations, which is a function of the values of $\partial U/\partial n$ on the cooling duct boundary. 

The coupled system is discretized using a finite elements scheme for the conduction part~\eqref{eq:conduction_differential}. 
This is done, because it ensures we can choose the mesh for the conduction part in a way that it is identical with the mesh used for the structural mechanics simulation described in Section \ref{ShapeOpStruc}, that uses the calculated temperatures as an input. 
For the cooling duct part \eqref{eq:cooling_differential}, 
we use a finite volume scheme, as that makes it easier to have energy conservation across the boundary and provides a clear mapping of the PDE boundary to cooling channel elements.
The resulting discretized system is then solved by solving each subsystem and updating the boundary condition respectively the right hand side of the other system, alternating between the two subsystems until the solutions of two consecutive iterations differ by a sufficiently small margin. 
This back-and-forth iteration is reminiscent of a Gauß-Seidel iteration scheme, or more general, a fixed-point iteration. 

\begin{figure}[ht]
    \centering
    \includegraphics[width=0.49\textwidth]{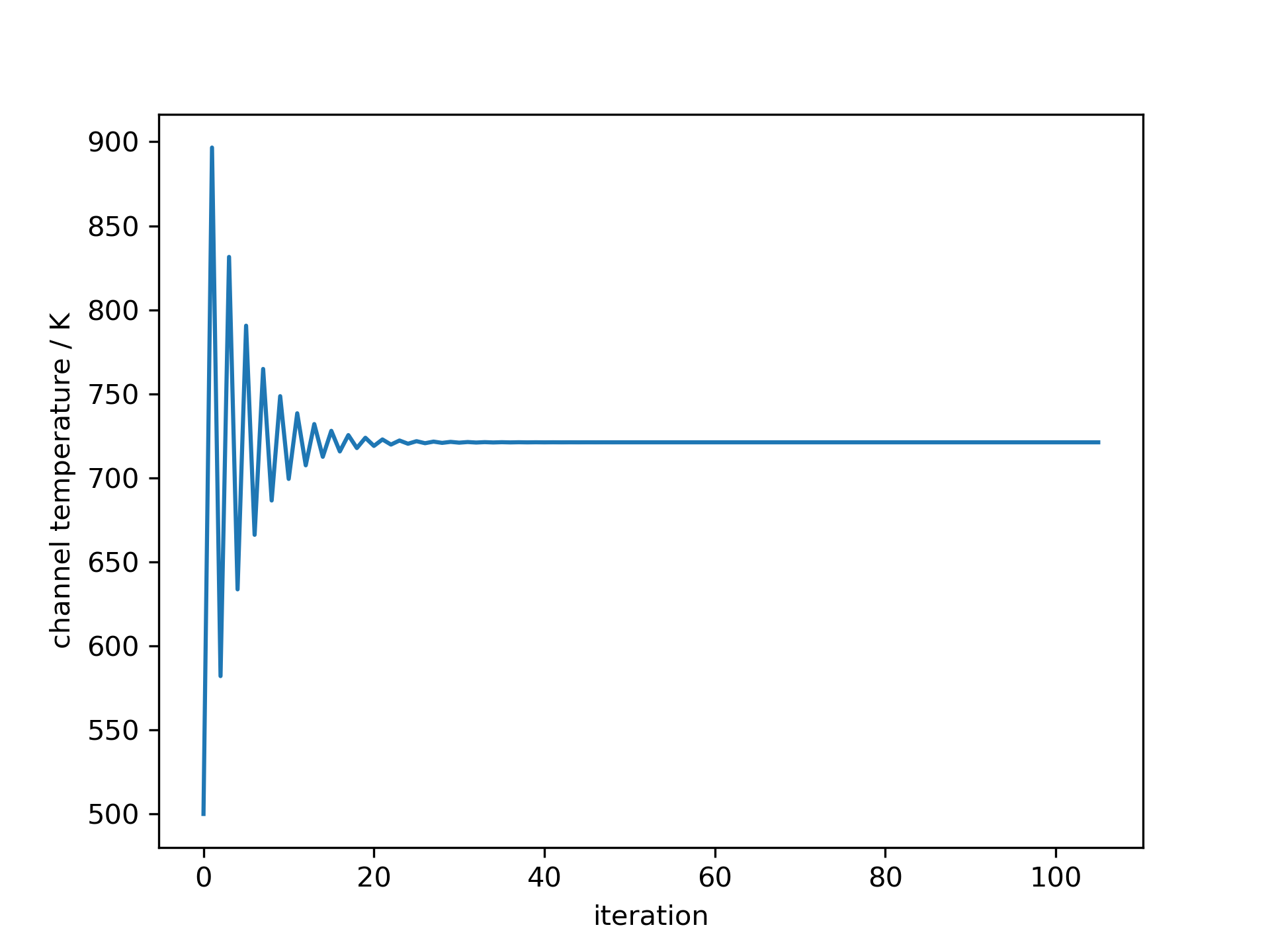}
    % \hspace{1em}
    \includegraphics[width=0.49\textwidth]{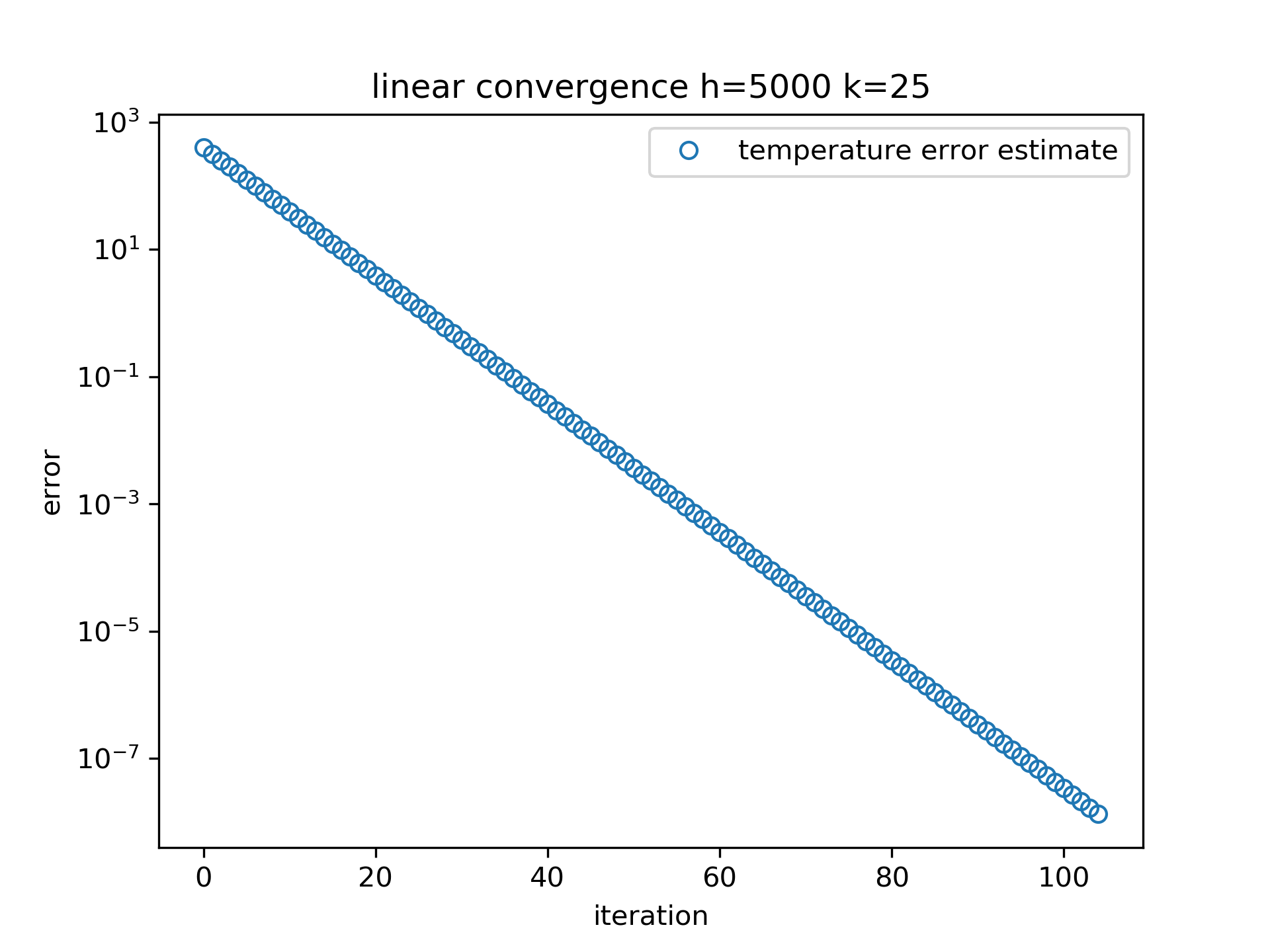}
    \caption{Behavior of the cooling duct outlet temperature (left) and error-estimate (right) for a converging set of parameters ($h_{int} = h_{ext} = 5000, k=25$)}
    \label{fig:cooling_channel_converging_solution}
\end{figure}

\begin{figure}[ht]
    \centering
    \includegraphics[width=0.49\textwidth]{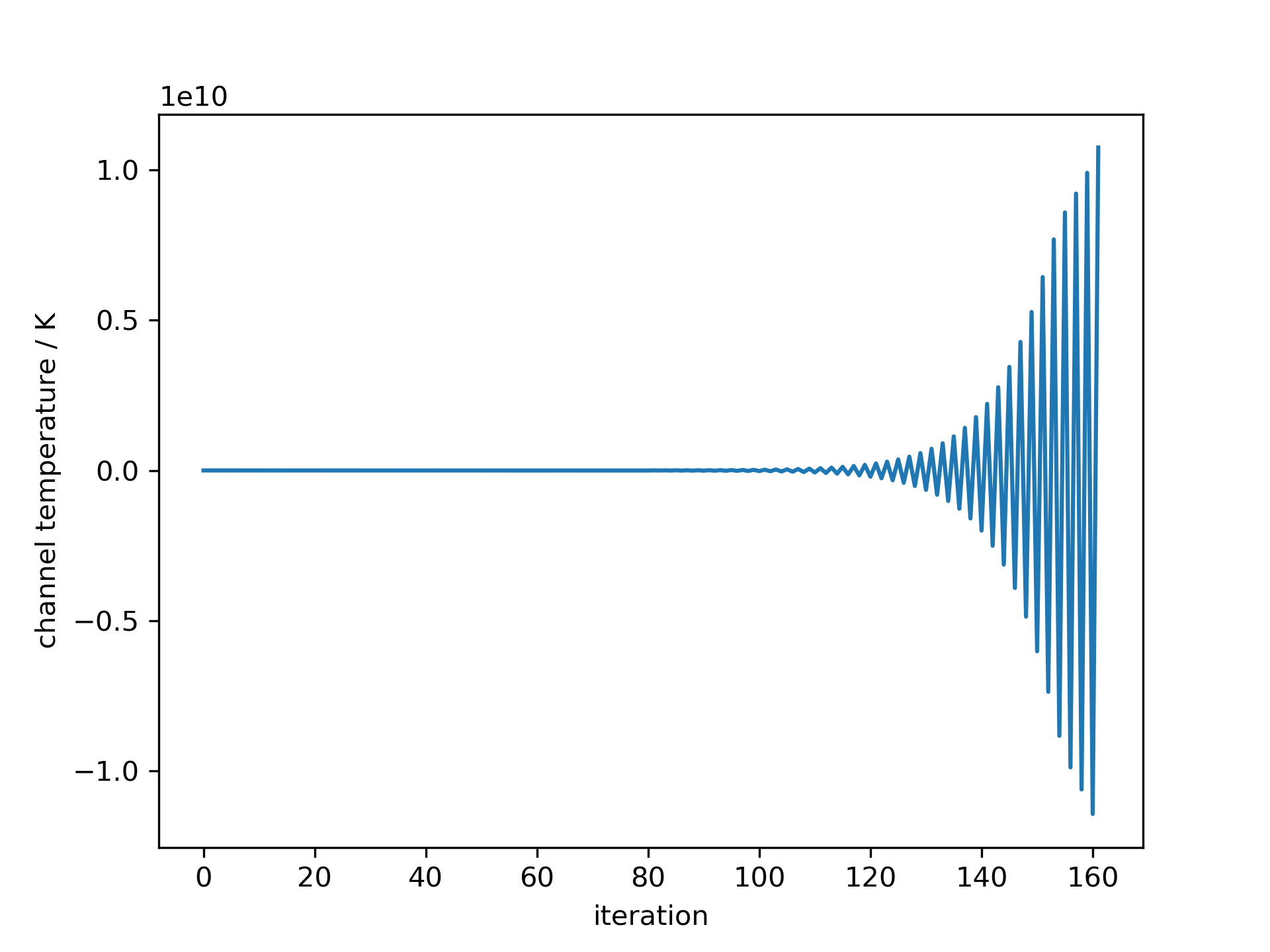}
    % \hspace{1em}
    \includegraphics[width=0.49\textwidth]{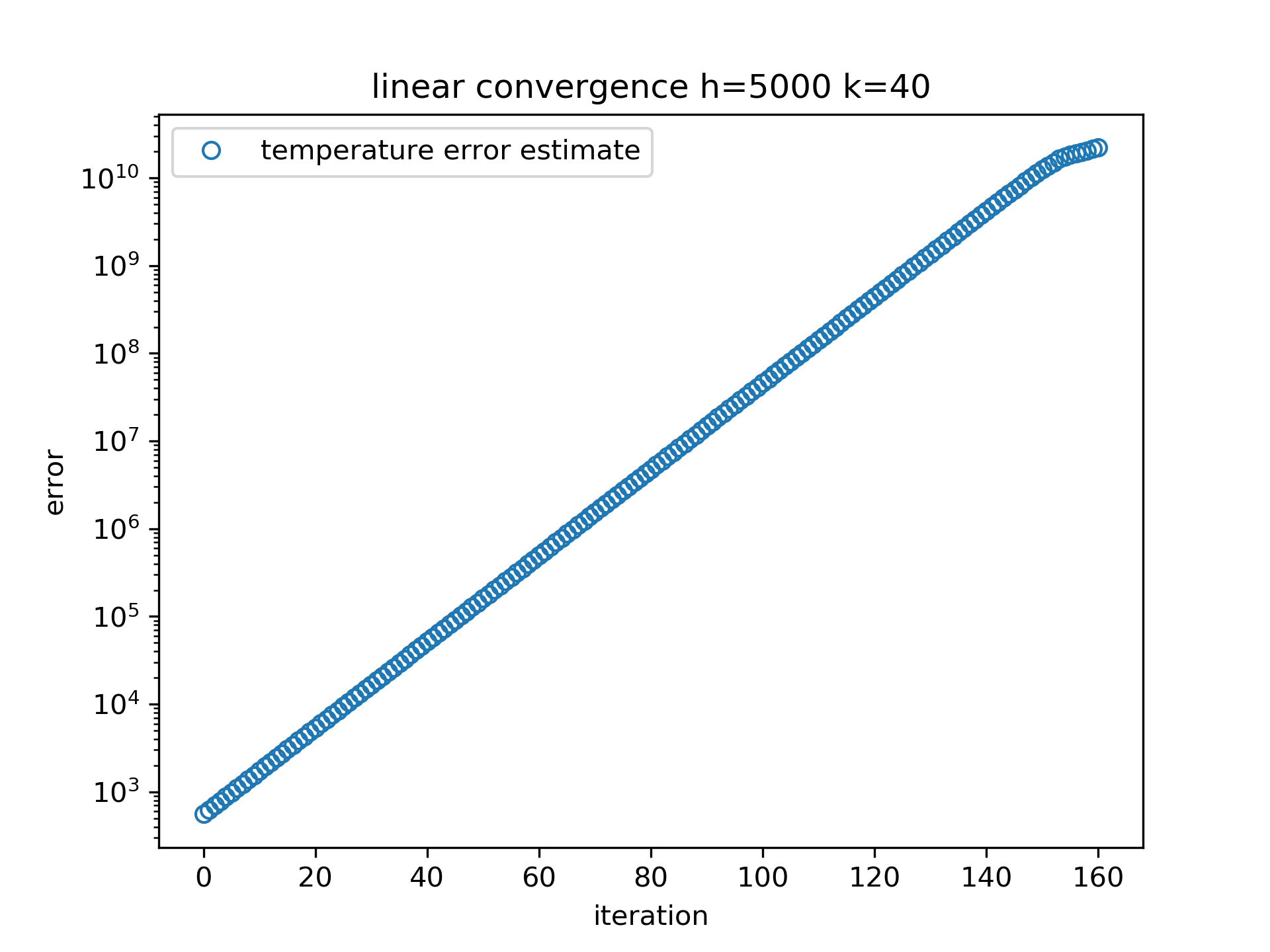}
    \caption{Behavior of the cooling duct outlet temperature (left) and error-estimate (right) for a diverging set of parameters ($h_{int} = h_{ext} = 5000, k=40$)}
    \label{fig:cooling_channel_diverging_solution}
\end{figure}

Numerical tests have shown that this iterative solution indeed exhibits linear convergence as seen in Fig.~\ref{fig:cooling_channel_converging_solution}, with the solution behaving like a dampened oscillation approaching the "correct" solution. 
These numerical tests also indicated that the convergence is not unconditional, but depends on the parameter values chosen for the system, especially the thermal conductivity $k$ and the heat transfer coefficients $h$. 
High values of $h$ lead to divergence and turn the aforementioned dampened oscillation into one with an exponentially increasing amplitude as seen in Fig.~\ref{fig:cooling_channel_diverging_solution}. 

%%%%%%%%%%%%%%%%%%%%%%%%%%%%%%%%%%%% WP3 (PI Gottschalk)
\subsection{Probabilistic Objective Functionals for Material Failure}\label{sec:23}
Since the pioneering work of Weibull~\cite{weibull1939}, the probabilistic modelling of material failure has been an established field of material science, see about \cite{batendorf1974}. 
Applications to the \textit{Low Cycle Fatigue} (LCF) damage mechanism can be found in \cite{olschewski1997, fedelich1998, sornette1992}.
In these studies, crack formation is modelled by percolation of intra-granular cracks or by kinetic theory for the combination of cracks.
The mathematical literature mainly contains generic volume or surface target functions without direct material reference. 
In numerical studies, global compliance is usually chosen as the objective functional, which also does not establish a direct relationship to material failure, see e.g.\ \cite{conti2011}.

The objective functional used in GivEn for the probability of failure originates \cite{gottschalk2014, gottschalk2014b, MaedeSchmitz_2018} see also \cite{gottschalk2015} for multi-scale modeling. 
A connection between probabilistic functional objectives of materials science and the mathematical discipline of shape optimization is produced for the first time in \cite{gottschalk2014}, 
see also \cite{bittner2016, bolten2015, Bolten2019}.

The aim in this sub-area is the probabilistic modelling of material damage mechanisms and the calculation of form derivatives and form Hessian operators for the failure probabilities of thermal and mechanically highly stressed turbine blades. 

The physical cause of the LCF mechanism in the foreground is the sliding of crystal dislocations along lattice planes with maximum shear stress and is therefore dependent on the random crystal orientation. For this reason, so-called intrusions and extrusions occur at the material surface, which eventually lead to crack formation 
\cite{baeker08, RadajVormwald_2007}.
In an effective approach, this scattering of material properties can be empirically investigated within the framework of reliability statistics.

In the deterministic approach prevalent in mechanical engineering, life expectancy curves are used to determines the service life at each point of the blade surface. 
The shortest of these times is to failure over all points is then, under consideration of safety discounts, converted to the permitted safe operating time of the gas turbine. 
The minimum formation inherent in this process means that the target functions cannot be differentiated. 
Probabilistic target functions, on the other hand, can be defined according to the form and continuously adjusted. 

In particular, the stability of discretization schemes must be examined in both with regard to geometric approximation of the forms as well as the solutions. 
The background
is that $H^1$ solutions are insufficient for a finite probabilistic target functional, especially if
notch support is also considered \cite{maede2017}. 
Suffice of this must be used a $W^{k,p}$ solution and approximation theory \cite{Ciarlet88}. 

Next, the calculation of the shape gradients and shape Hesse operators of the functionals essentially follows \cite{Sokolowski92}, with open questions about the existence and properties of shape gradients for the surface- and stress-driven damage mechanism LCF still to be clarified. This program has been started within the GivEn research initiative, 
cf.\ \cite{bittner2020shape}. 
Analogous to \cite{gottschalk2014}, the solution strategy is based on a uniform
regularity theory for systems of elliptic PDEs, 
cf.\ \cite{agmon1959, agmon1964, Ciarlet88}.
In particular, the mathematical status of the continuously-adjusted equation deserves further attention,
as this has a high regularity loss for surface-driven LCF. 

In the following we present a hierarchy of probabilistic failure models that give rise to objective functionals related to reliability. We start with the simple Weibull model, proceed with a probabilistic model for LCF proposed by \cite{SchmitzSeibel_2013,MaedeSchmitz_2018} and then give an outlook on the multi-scale modeling of the scatter in probabilistic LCF, see \cite{engel2019probabilistic}.

\subsubsection{The Weibull Model via Poisson Point Processes}
Technical ceramic has multiple properties such as heat or wear resistance that make them a widely used industrial material. 
Different to other industrial material, the physical properties of ceramic materials highly depend on the manufacturing process. 
What determines the failure properties the most, are small inclusions that stem from the sintering process. 
These make ceramic a brittle material, leading to a somewhat high possibility of failure of the component under tensile load often before the ultimate tensile strength is reached \cite{baeker08}. 

When applying tensile load, these inclusions may become the initial point of a crack, developing into a rupture if a certain length of the radius of the crack is exceeded at a given level of tensile stress. 
Therefore, the probability of failure under a given tensile load is the probability that a crack of critical length occurs. Or to phrase it differently, the survival probability in this case is the probability that exactly zero of these critical cracks occur. Thus, for a given domain $\Omega\subset\mathbb{R}^d$, $d=2,3$ with a suitable counting measure $\mathcal{N}$ \cite{kallenberg}, we can express the failure probability in the following way,
\begin{align}
    \text{PoF}(\Omega) =1- \mathbb{P}\bigl(\mathcal{N}( A_c(\Omega )) =0\bigr),
\end{align}
where $A_c(\Omega ))$ is the set of critical cracks.
The probability, that one of the inclusions grow into a critical crack, mainly depends on the local stress tensor $\sigma_n(u)$, which itself is determined by a displacement field $u\in H^1(\Omega,\mathbb{R}^d)$, that is the solution of a linear elasticity equation.
As there is no other indication, it is feasible to assume that the location, size and orientation of the initial inclusions are independent of each other and uniformly randomly distributed. Under these assumptions, it follows that the counting measure $\mathcal{N}(\Omega)$ is a Poisson point process (PPP). 
Taking further material laws into account it follows that \cite{bolten2015}
 \begin{align} \label{eqa:PSurvival}
 \text{PoF}(\Omega |u)=1-\mathbb{P}(\mathcal{N}(A_c(\Omega ,u))=0)
 =1-\exp\{-\nu (A_c(\Omega, u))\},
 \end{align}
 with the intensity measure of the PPP
 \begin{align}
     \nu (A_c(\Omega, u))
     =\frac{\Gamma(\frac{d}{2})}{2\pi^{\frac{d}{2}}}\int\limits_{\Omega}\int\limits_{S^{d-1}}\int\limits_{a_c}^{\infty}d\nu_a(a)\,dndx.
\end{align}
With some reformulations we find our objective functional 
of Weibull type
\begin{align}\label{eqa:ObFun}
J_1(\Omega ,u):=\nu (A_c(\Omega,u))=\frac{\Gamma(\frac{d}{2})}{2\pi^{\frac{d}{2}}}\int\limits_{\Omega}\int\limits_{S^{d-1}}\left(\frac{\sigma_n}{\sigma_0}\right)^m\,dndx.
\end{align}
This functional~\eqref{eqa:ObFun} will be one of the objective functionals in the following (multiobjective) gradient based shape optimization.

%%%%%%%%%%%%%%%%%%%%%%%%%%%%%%
\subsubsection{Probabilistic Models for LCF} 
Material parameters relevant for fatigue design, like the HCF fatigue resistance were considered as a random variable for a long time \cite{Siebert_1956,Murakami_1989,RadajVormwald_2007} and distributions and their sensitivities were even recorded in general design practice standards \cite{Monnot_1988,hertel12,FKMnlin_2015,FKMnlin_2018}. 
The existence of flaws, such as crystal dislocations, non-metallic inclusions or voids in every material has early lead to the discovery of the statistical size effect 
\cite{Phillipp1942,Markovin1944,Massonet1956,Hempel1957}. 
Within the last decade, a local probabilistic model for LCF based on the Poisson point process was 
developed by Schmitz \textit{et al.} \cite{SchmitzSeibel_2013,SchmitzRollmann_2013} for predicting the statistical size effect in any structural mechanics FEA model. 
It approximates the material LCF life statistics with a Weibull distribution which allows developing a closed form integral solution for the distribution scale $\eta$ (see equation~\eqref{eq:OptRelia_ScaleSurfInt}). 
Recently, M\"ade \textit{et al.} \cite{MaedeSchmitz_2018} presented a validation study of the combined size and stress gradient effect modeling approach within the framework of Schmitz \textit{et al.} \cite{SchmitzSeibel_2013}.
If stress gradients are present in components, they can have an increased LCF life. The benefit is proportional to the stress gradient but also material dependent  \cite{ErikssonSimonsson_2018,SiebelMeuth_1949,Siebel_1955,Wundt_1972}.
%Since those are typically related to localized stress concentrations, a beneficial size effect simultaneously shifts the LCF life distribution scale parameter to higher values. 
A stress gradient support factor $n_\chi=n_\chi(\chi\mid\boldsymbol\vartheta)$, a functional of the normalized stress gradient $\chi$ and material-specific parameters $\boldsymbol\vartheta$ \cite{Siebel_1955}, is introduced to quantify the effect.
While the size effect as described by the surface integral~\eqref{eq:OptRelia_ScaleSurfInt} causes an actual delay or acceleration in fatigue crack initiation, researchers share the interpretation that stress gradient support effects in LCF root back to retarded propagation of meso-scale cracks in the decreasing stress field \cite{RadajVormwald_2007,LazzarinTovo_1997,MansonHirschberg_1967,KontermannAlmstedt_2016,Dowling_1979,SakaneOhnami_1986}. 
%In experiments, such LCF cracks are therefore sometimes detected after many more cycles than predicted with the statistical size effect \cite{MaedeSchmitz_2018,SchmitzPhD_2014}. That is why $n_\chi$ is much dependent on the crack detection criterion \cite{KontermannAlmstedt_2016,FFEMS_MaedeKumar_2019,Mughrabi_2015}. 
In order to compare the stress gradient effect for different materials, a common detectable, ``technical'' crack size must be defined. 
%LCF notch specimen test data can then be used for calibrating the stress gradient support factor. 
Since the stress gradient $\chi(\mathbf{x})$ is, like the stress field, a local property, it was integrated into the calculation of the local deterministic life $N_\mathrm{det}(\mathbf{x})$ with the 
Coffin-Manson-Basquin model:
\begin{align}
  \frac{\varepsilon_a(\mathbf{x})}{n_\chi\left(\chi(\mathbf{x})\mid\mathbf{\vartheta}\right)}
  =\frac{\sigma^\prime_f}{E}\cdot\left(2N_\mathrm{det}(\mathbf{x})\right)^b+\varepsilon^\prime_f\cdot\left(2N_\mathrm{det}(\mathbf{x})\right)^c \label{eq:ProbMatSc_CMB_NoS}.
\end{align}
Here, the stress is computed with the aid of the linear elasticity equation, which this time is not a technical tool for smoothing gradients as in \eqref{eq:elas}, but represents the physical state, namely
\begin{equation*}\label{eq:el-eq}
\begin{tabular}{ l l }
$\nabla \cdot \sigma(u) + f = 0$ & in $\Omega$ \\
$\sigma(u) = \lambda(\nabla \cdot u) I + \mu (\nabla u + \nabla u^\top)$ & 
in $\Omega$ \\
$u = 0$ & on $\partial\Omega_D$ \\
$\sigma(u) \cdot n = g$ & on $\partial\Omega_N$. \\
\end{tabular}
\end{equation*}
Here, $\Omega$ represents the component, $\lambda>0$ and $\mu>0$ are Lam\'e
coefficients and $u:\Omega\to\mathbb{R}^3$ is the displacement field on $\Omega$ obtained as a reaction to the volume forces $f$ and the surface loads $g$. 
We connect the topic of optimal probabilistic 
reliability to shape optimization elasticity PDE as state equation and classify Poisson point process models according to their singularity \cite{bittner2020shape}.
Following \cite{SchmitzSeibel_2013,MaedeSchmitz_2018}, we obtain for the probability of failure at a number of use cycles $n$
\begin{equation}
    \text{PoF}(\Omega,n)=1-e^{-n^m J_R(\Omega,u)}
\end{equation}
The functional $J(\Omega,u)$ that is arising out of this  framework is given by:
\begin{equation}
	J_R(\Omega, u_\Omega):=
	\int_{\partial\Omega} \left( 
	\frac{1}{N_\text{det} (\nabla
		u_\Omega(x),\nabla^2
		u_\Omega(x))}
	\right)^m \,dA. \label{eq:OptRelia_ScaleSurfInt}
\end{equation}
$N_\text{det}$ denotes the deterministic numbers of
life cycles at each point of the surface of the component and 
$m$ is the Weibull shape parameter. 

\begin{figure}[t]
	\centering
	\includegraphics[width=0.95\textwidth]{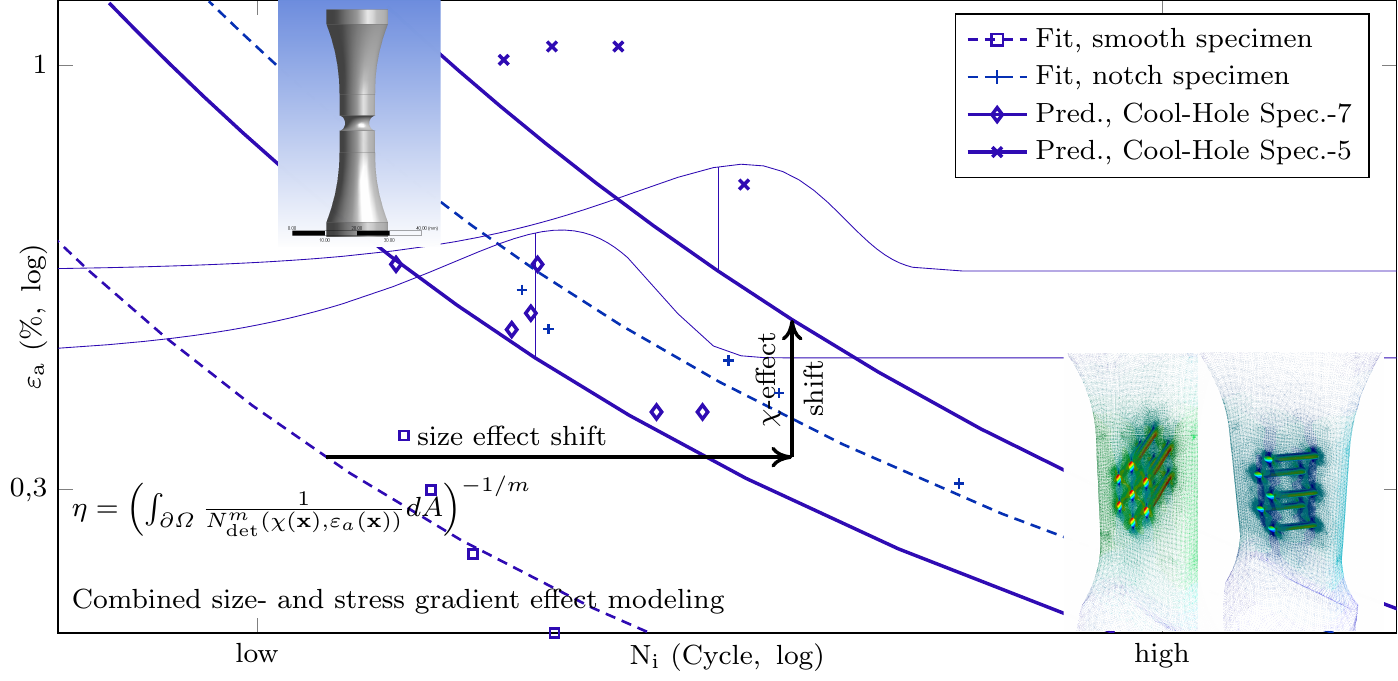}
	\caption{Strain W\"ohler plot of LCF test data, calibrated (dashed) and predicted (solid) median curves for smooth ($\Box$), notch ($+$) and cooling hole specimens ($\diamondsuit,\,\times$). 
	All W\"ohler curves are interpolated median values of the Weibull LCF distributions exemplary indicated with the thin density function plot.}
	\label{fig:ProbMatSc_ComEffValidat}
\end{figure}

M\"ade \textit{et al.} have calibrated the material parameters $\mathbf{\vartheta}$, $E$, $\sigma^\prime_f$, $\varepsilon^\prime_f$, $b$, $c$ as well as the Weibull shape parameter $m$ with the Maximum-Likelihood method simultaneously using smooth and notch specimen data simultaneously \cite{maede2017,MaedeSchmitz_2018}. The resulting model was able to predict the LCF life distribution for certain component-similar specimens (see Fig.~\ref{fig:ProbMatSc_ComEffValidat}\footnote{Reprinted from Comp. Mat. Sci., 142, (2018) pp. 377–388, M\"ade et al., Combined notch and size effect modeling in a local probabilistic approach for LCF, Copyright (2017), with permission from Elsevier}).

In the following, we  apply this model as cost functional in order to optimize the component $\Omega$ w.r.t.\ reliability.

%%%%%%%%%%%%%%%%%%%%%%%%%%%%%%%%%%%%%%%%%%%%%%%%%
\subsubsection{ Multi-Scale Modeling of Probabilistic LCF}

While the Weibull-based approach from the previous subsection allows a closed-form solution and therefore fast risk assessment computation times, 
the microstructural mechanisms of LCF suggest a different distribution shape \cite{MochPhD_2018}. 
Since this is not yet assessable by LCF experiments in a satisfying way, Engel \textit{et al.} have used numerical simulations of probabilistic Schmid factors to create an LCF model considering the grain orientation distribution and material stiffness anisotropy in cylindrical Ni-base superalloy specimens \cite{engel2019probabilistic}.

\begin{figure}[H]
	\centering
	\includegraphics[width=0.95\textwidth]{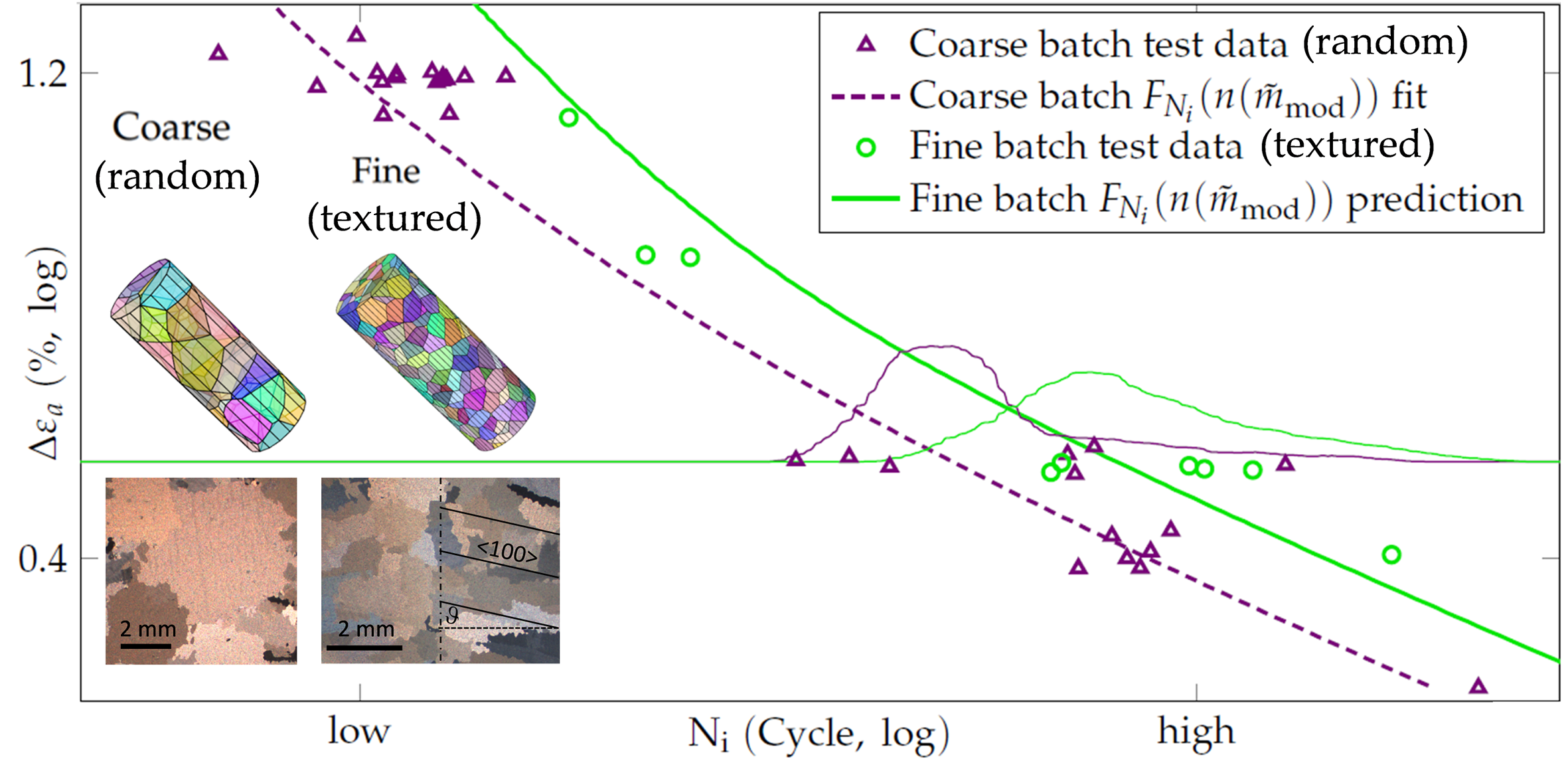}
	\caption{Strain W\"ohler plot of LCF test data, calibrated (dashed) and predicted (solid) median curves for specimens with isotropic and preferential grain orientation distribution (coarse and fine grain) by Engel \textit{et al.} \cite{engel2019probabilistic}. Calibration and prediction was carried out using the Schmid factor based LCF life distribution. 
	The underlying Schmid factor distribution was derived from polycrystalline FEA simulations which considered the lattice anisotropy and orientation distribution in both specimen types.}
	\label{fig:ProbMatSc_SchmidLCFPred}
\end{figure}

Polycrystalline FEA models were developed to investigate the influence of local multiaxial stress states a result of as grain interaction on the resulting shear stress in the slip systems. Besides isotropic orientation distributions also the case of a preferential orientation distribution was analysed.
From the FE analyses, a new Schmid factor distribution, defined by the quotient of $\mathrm{max}(\tau_{rss})$ maximum resolved shear stress at the slip systems and von Mises stress $\sigma_{vM}$, was derived as a probabilistic damage parameter at each node of the model. 
Qualitatively as well as quantitatively, they differ largely from the single grain Schmid factor 
distribution of Moch~\cite{MochPhD_2018} and also from the maximum Schmid factor distribution of grain ensembles presented by Gottschalk \textit{et al.} \cite{GottschalkSchmitz_2015}. 
Experimental LCF data of two different batches presented by 
Engel~\cite{EngelPhD_2018,EngelLCF8_2017} showed different LCF resistances and microstructural analyses revealed a preferential grain orientation in the specimens which withstood more cycles (see Fig.~\ref{fig:ProbMatSc_SchmidLCFPred}). 
By combining the Schmid factor based LCF life model of 
Moch~\cite{MochPhD_2018} and the Schmid factor distribution generated by FEA, Engel \textit{et al.} were able to predict just that LCF life difference \cite{engel2019probabilistic}. 
Ultimately, it was found that the microstructure-based lifing model is able to predict LCF lives with higher accuracy than the Weibull approach by considering the grain orientation and their impact on the distributions of Young’s moduli and maximum resolved shear stresses. 
However, the application is computationally demanding and its extension to arbitrary components still has to be validated.
% Die Autoren halten das Copyright (CC-BY Lizenz)

%%%%%%%%%%%%%%%%%%%%%%%%%%%%%%%%%%%%%%%%%% WP4 (PIs Bolten)
\subsection{Shape Optimization for Probabilistic Structure Mechanics}\label{ShapeOpStruc}
This sub-area deals with two kinds of failure mechanisms, failure of brittle material under tensile load and low cycle fatigue (LCF). 
As explained in the preceding section, the probability of failure for both failure mechanisms can be expressed as local integral of the volume (ceramics) or surface (LCF) over a non-linear function that contains derivatives of the state  $u_\Omega$ subject to the 
elasticity equation~\eqref{eq:el-eq}.

%Hence, the failure probability of these two kinds have to be modeled in different ways.
%As a brittle material, we consider ceramic material. 

Since problems in shape optimization generally do not result in a closed solution, the numerical solution plays a major role, e.g.\ in \cite{Haslinger03}. 
Typically, the PDEs occurring as a constraint are discretized using the finite element method. 
Here, a stable and accurate mesh representation of $\Omega$ is needed for stable numerical results. 
Deciding for a mesh always means balancing the need for accuracy of the representation on the one hand and on the other hand the time to solution.
Especially in applications such as shape optimization, 
where usually hundreds to thousands of iterations and thus changes in the geometry and mesh are needed to find a converged solution, the meshing in each iteration often becomes a bottle neck in terms of computational cost.
Recent research therefore aims to find methods to move the grid points of a given representation in a stability preserving way, rather than to perform a re-meshing. 
These approaches result in unstructured grids. 
To exploit the means of high performance computing however, 
structured grids are way more desirable.
This let us to consider an approach first developed in \cite{hackbusch1996, hackbusch1997a, hackbusch1997b}. 
For demonstration purposes, the technique is described in two dimensions but easily extends to three dimensions.

\begin{figure}[t]
	\centering
		\subfloat{\includegraphics[width=0.3\linewidth]{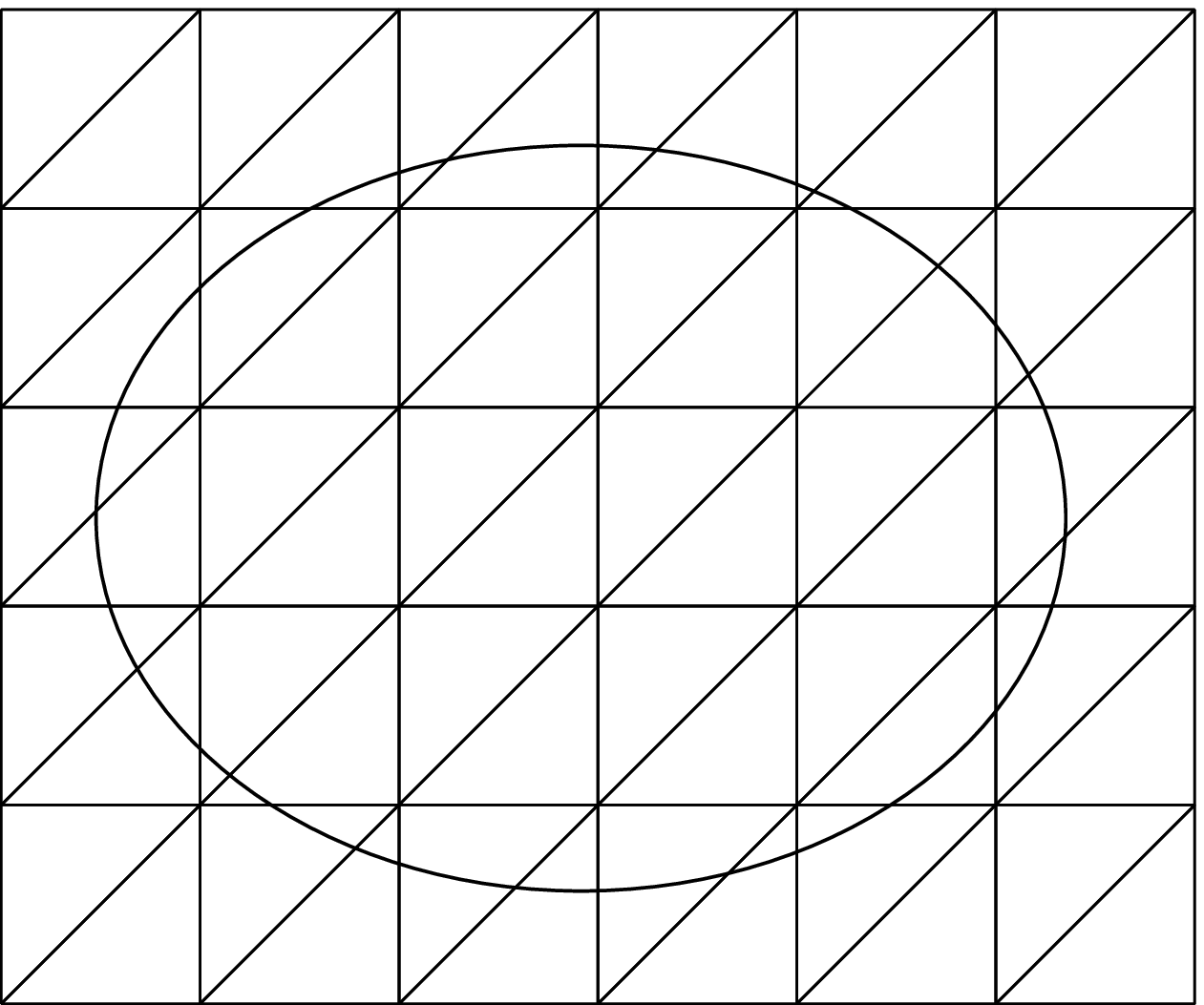}}
		\hfill
		\subfloat{\includegraphics[width=0.3\linewidth]{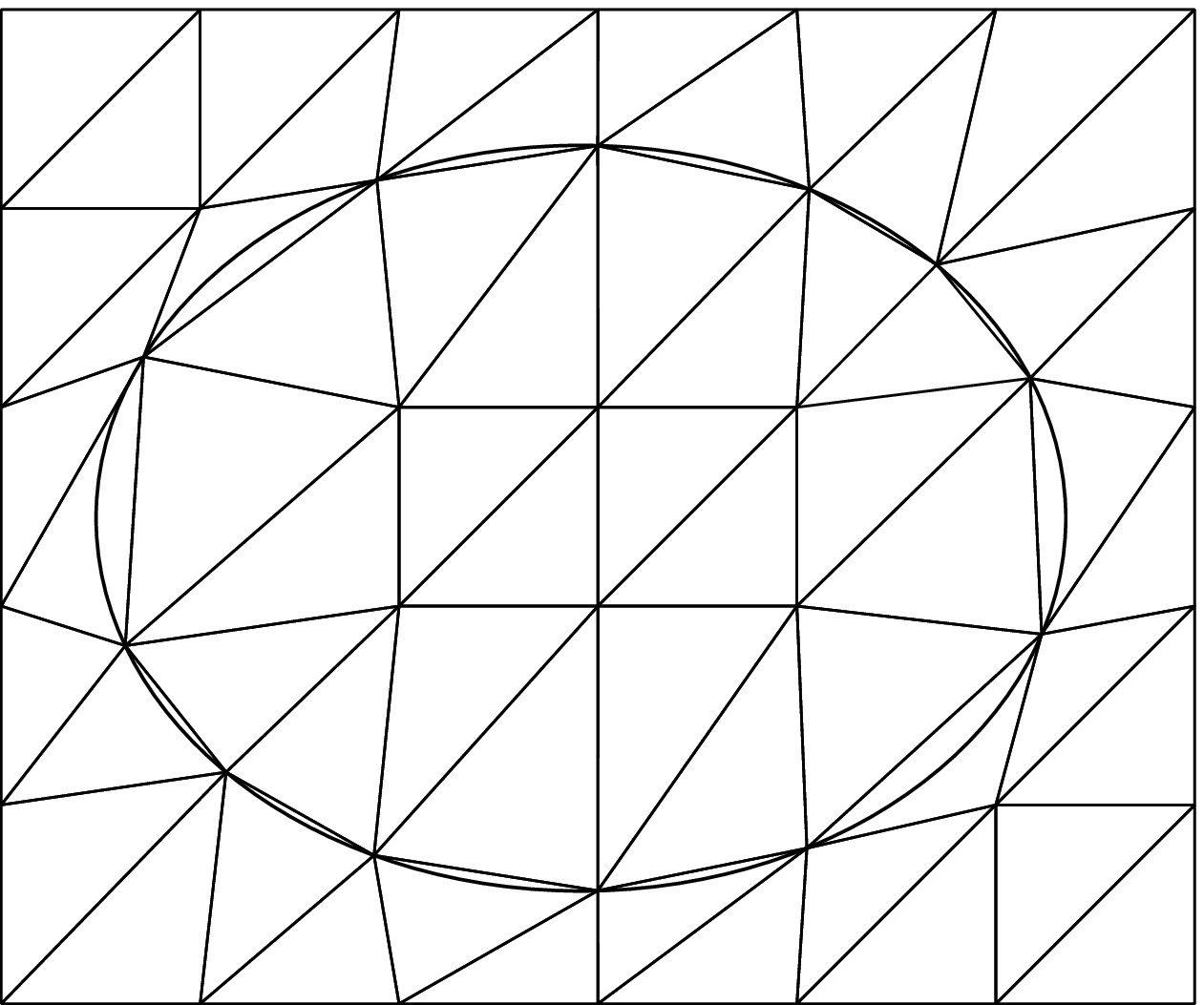}}
		\hfill
		\subfloat{\includegraphics[width=0.3\linewidth]{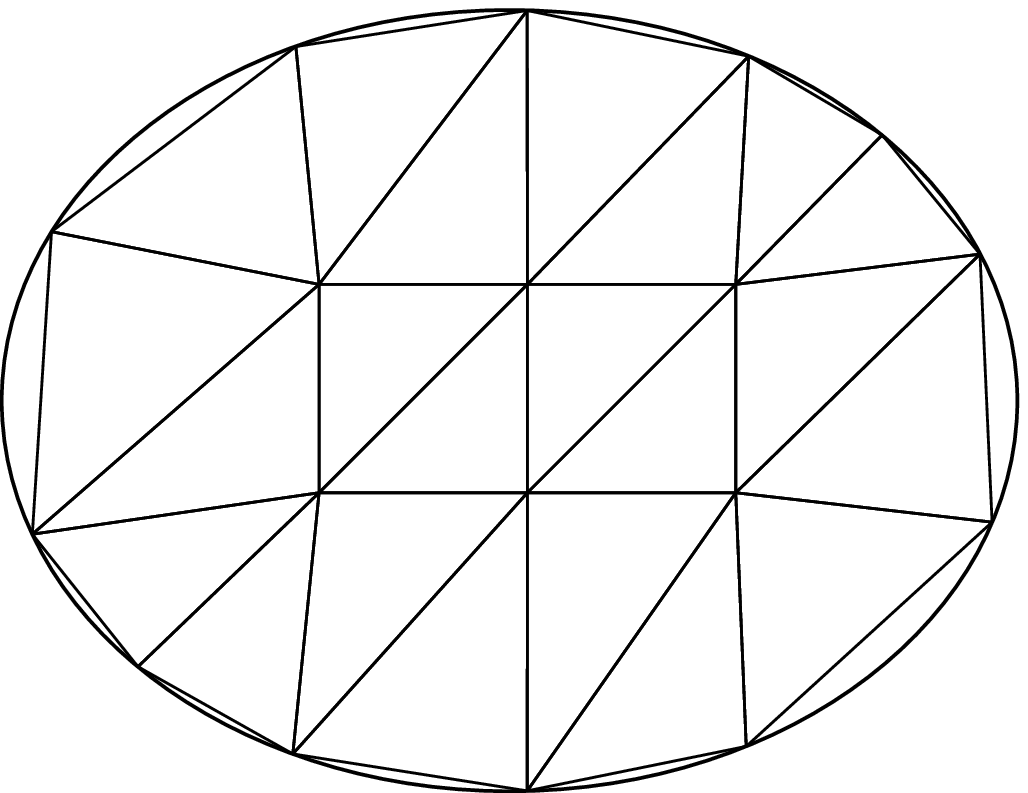}}
	\caption{Visualization of adaption of the grid\protect\footnotemark}
	% use \protect in front of \footnotemark in captions
	\label{fig:cfe}
\end{figure}

\footnotetext{Reprinted from Progress in Industrial Mathematics at ECMI~2018, pp. 515-520, Bolten et al., Using Composite Finite Elements for Shape Optimization with a Stochastic Objective Functional, Copyright: Springer Nature Switzerland AG 2019}
We consider a rectangular domain $\tilde{\Omega}$ which is discretized by a regular triangular grid. 
We assume that all admissible shapes that occur during the optimization process lie in this domain. 
The regular grid on $\tilde{\Omega}$ is denoted by $\tilde{T}$, the number of elements by $N^{el}$ and the number of nodes by $N^{no}$.
In a second step, the boundary $\delta\Omega_0$ of the shape to be optimized $\Omega_0$ is superimposed onto the grid, see Fig.~\ref{fig:cfe}a. 
The regular grid is then adapted to the boundary by moving the closest nodes to the intersections of the grid and the boundary, see Fig.~\ref{fig:cfe}b. 
The adapted grid is denoted by $\tilde{T}_0$.
As the nodes are moved only, the connectivity of $\tilde{T}_0$ is the same as before. 
During the adaption process, the cells of $\tilde{T}_0$ are assigned a status as cells lying inside or outside of the component. 
The computations are only performed on those cells that are inside the component. When updating the grid according to the new shape of the domain $\Omega_1$ and so forth, the adaption process starts with $\tilde{\Omega}$ again, while taking into account the information about the previous domain, the step length and search direction that led to the current domain. 
By this, only the nodes that lie in a certain neighborhood of the previous boundary have to be checked for adaption. 
This leads to a speed up in meshing compared to actual re-meshing techniques. 
Additionally, as the grid is otherwise regular and the connectivity is kept, the governing PDE only has to be changed in entries representing nodes on elements that have been changed, hence no full assembly is needed, which is an advantage to both re-meshing and mesh morphing techniques.

The objective functional~\eqref{eqa:ObFun} is discretized via finite elements. For reduction of the computational cost, the adjoint approach than leads to the derivative
\begin{align}\label{eq:grad}
 \frac{dJ_1\bigl(X,U(X)\bigr)}{dX} &= 
 \pp{J_1(X,U)}{X}
 +\Lambda\left[
 \pp{F(X)}{X}-\pp{ B(X)}{X}U\right]\\
 \label{eq:adjoint}
B^\top(X)\Lambda &= \pp{J_1(X,U)}{U}\\
\label{eq:elas}
B(X)U &= F(X),
\end{align}
with \eqref{eq:adjoint} being the adjoint equation giving the adjoint state and \eqref{eq:elas} is the discretized linear elasticity equation, giving the discrete displacement $U$. $X$ represents the discretized domain $\Omega$.

\begin{figure}[t]
	\centering
		\subfloat[Standard gradient]{\includegraphics[width=0.5\linewidth]{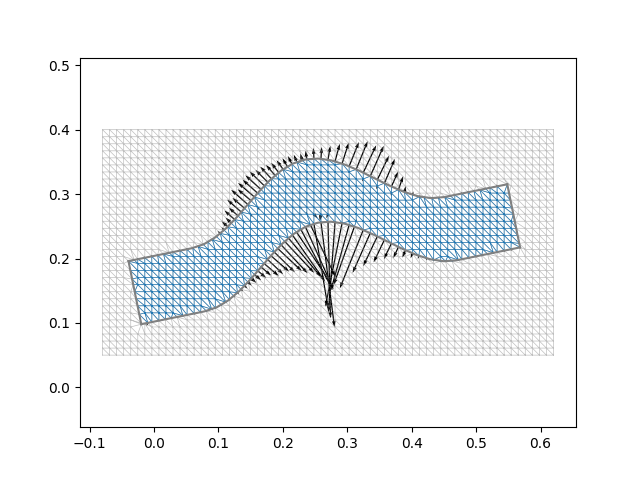}}
		\hfill
		\subfloat[Smoothed gradient]{\includegraphics[width=0.5\linewidth]{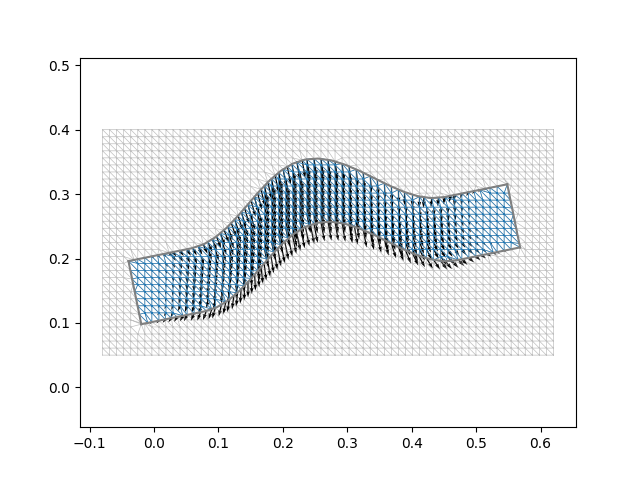}}
	\caption{Gradient in standard scalar product and smoothed gradient; $N_x=64, N_y=32$}
	\label{fig:ceramic-gradients}
\end{figure}

With (\ref{eq:elas}) and (\ref{eq:adjoint}) the derivative (\ref{eq:grad})  is calculated on the structured mesh as visualized in Fig. \ref{fig:ceramic-gradients}a. For the optimization, more closely described in the following section \ref{sec:25}, the gradient is smoothed using a Dirchlet-to-Neumann map \cite{Schmidt09}
(see Fig. \ref{fig:ceramic-gradients}b). This provides the shape gradients needed for further gradient based optimization steps in the following section.

%%%%%%%%%%%%%%%%%%%%%%%%%%%%%%%%%%%%%%%%%% WP5 (PIs Klamroth and Stiglmayr)
\subsection{Multiobjective Optimization}\label{sec:25}
The engineering design of complex systems like gas turbines often requires the consideration of multiple aspects and goals. 
Indeed, the optimization of the reliability of a structure usually comes at the cost of a higher volume and, hence, a higher production cost. 
Other relevant optimization criteria are, for example, the minimal buckling load of a structure or its minimal natural frequency \cite{Haslinger03}. 
In this section, we consider both the mechanical integrity and the cost of a ceramic component in a biobjective PDE constrained shape optimization problem. 
Further objective functions can (and should) be added to the model depending on the application at hand.  
Towards this end, we model the mechanical integrity $J_1(\Omega,u)$ of a component $\Omega$ as described in Section~\ref{ShapeOpStruc}, see \eqref{eqa:ObFun}, while the cost $J_2(\Omega)$ is assumed to be directly proportional to the volume of the component, i.e.,  $J_2(\Omega)=\int_{\Omega} dx$. 

Multiobjective shape optimization including mechanical integrity as one objective is widely considered, see, e.g., \cite{chir:mult:2018} for a recent example. Most of these works neither consider probabilistic effects nor use gradient information. The formulation introduced in Section~\ref{ShapeOpStruc} overcomes these shortcomings. 
It was first integrated in a biobjective model in \cite{Doganay2019}, where two alternative gradient-based optimization approaches are presented. We review this approach and present new numerical results based on structured grids and advanced regularization.

%%%%%%%%%%%%%%%%%%%%%%%%%%%%%%%%%%%%%%%%%%%%%%%%%%
\subsubsection{Pareto Optimality} 
\label{pareto_optimal}
Multiobjective optimization asks for the simultaneous minimization of $p$ conflicting  objective functions $J_1, \dots, J_p$, with $p\ge2$. 
We denote by $J(\Omega)=(J_1(\Omega),\dots,J_p(\Omega))$ the outcome vector of a feasible solution  $\Omega\in\Oad$ (i.e., an admissible shape).
Since in general the optimal solutions of the objectives $J_1, \dots, J_p$ do not coincide, a multidimensional concept of optimality is required. 
The so-called Pareto optimality is based on the component-wise order \cite{Ehrgott2005}: 
A solution $\Omega \in \Oad$ is \emph{Pareto optimal} or \emph{efficient}, if there is no other solution $\Omega^{\prime} \in \Oad$ such that $J(\Omega^{\prime}) \leqslant J(\Omega)$, i.e., $J_i(\Omega^{\prime}) \le J_i(\Omega) \text{ for } i=1,\dots,p \text{ and } J(\Omega^{\prime}) \ne J(\Omega)$. 
In other words, a solution is efficient if it can not be improved in one objective $J_i$ without deterioration in one other objective function $J_k$.

\subsubsection{Foundations for Multi-Physics Multi-Criteria Shape Optimization} 
The existence of Pareto fronts for the multi-criteria case are  considered in a simplified analytical model replacing the RANS equations by potential theory with boundary layer losses. 
Pareto fronts can be replaced by scalarizations
using the techniques from \cite{gottschalk2014, bittner2016} or \cite{fuji1988, bolten2015}. 
Their continuous course is investigated by variation of the scalarization and the associated optimality conditions. 
The convergence of the discretized Pareto optimum solutions against the continuous Pareto optimum solution shall be studied according to the approach of \cite{Haslinger03}.

Optimizing the design of some component in terms of 
reliability, efficiency and pure performance includes the consideration of various physical systems that interact with each other. 
This leads naturally to a multicriteria shape optimization problem over a shape space $\mathcal{O}$ with requirements represented by cost functionals  $\mathbf{J} = (J_1, \dots, J_l)$:

\begin{equation*}
\begin{split}	
\begin{cases}
\text{Find } \Omega^* \in \mathcal{O} \text{ such that}\\
(\Omega^*,\boldsymbol{v}_{\Omega^*})\text{ is Pareto optimal w.r.t. }
\boldsymbol{J}.
\end{cases}
\end{split}
\end{equation*}

The class of cost functional we use to model 
the requirements on the component are connected
with physical state equations and 
arise from the probabilistic framework. They are described by

\begin{equation*}
		\begin{split}
		&J_{\text{vol}}(\Omega,\boldsymbol{v}) := \int_\Omega 
		\mathcal{F}_{\text{vol}}\hspace{.1em} 
		(x,\boldsymbol{v},\nabla \boldsymbol{v},\dots, \nabla^k \boldsymbol{v})\,dx, \\
		&J_{\text{sur}}(\Omega,\boldsymbol{v}) := \int_{\partial \Omega} 
		\mathcal{F}_{\text{sur}}\hspace{.1em}
		(x,\boldsymbol{v},\nabla \boldsymbol{v},\dots, \nabla^k \boldsymbol{v})\,dA.
	\end{split}
\end{equation*}
Due to the fact that the event of failure, e.g.\ the crack initiation process takes place on the surface of the component, the physical systems need to fulfill regularity conditions in order to be includable in this setting. 
We describe the possible designs of the component in this shape optimization problem by Hölder-continuous functions which give us the possibility to freely morph the shapes in various designs while remaining the premised regularity conditions. 
In this situation uniform regularity
estimates for the solutions of the physical system are needed in order to ensure the
existence of a solution to this design problem in terms of Pareto optimality.
The aim of this subproject is to translate a multi physical shape optimal design problem into the context of a well-posed multicriteria optimization problem. 

We couple internal and external PDEs in order to describe the various forces that are inflected on the
component. In this framework, using techniques based on pre-compactness of embedding between H\"older spaces of different index like in \cite{gottschalk2014,bittner2016}, we are able to show \cite{reese2020} the existence of Pareto optimal shapes in terms of
subsection \ref{pareto_optimal} which form a Pareto front, see also \cite{Doganay2019} for a related result.  We also prove the completeness of the Pareto front in the sense that the Pareto front coincides with the Pareto front of the closure of the feasible set (which is equivalent to the fact that every non-Pareto admissible shape is dominated by a Pareto optimal shape).

Further we investigated scalarization
techniques which transform the multi-objective optimization problem into a uni-variate problem.
In particular we considered the so-called achievement function and $\epsilon$-constraint methods
which depend, besides on the cost functional, on an additional scalarization parameter that represents 
the different weightings of the optimization targets, as e.g. reliability or efficiency. 
Hence, the shape space on which the optimization process takes place can also depend on this parameter and
with it the corresponding space of optimal shapes as well. 
Under suitable assumptions on the contuinuous dependency of the scalarization method on the scalarization parameter, 
we are able to show a continuous dependency 
of the optimal shapes spaces on the parameter as well. For details we refer to the forthcoming work \cite{reese2020}.

%%%%%%%%%%%%%%%%%%%%%%%%%%%%%%%%%%%%%%%
\subsubsection{Multiobjective Optimization Methods}
\label{sec:mop_methods}
Algorithmic approaches for multiobjective optimization problems can be associated with two common paradigms: scalarization methods and non-scalarization methods. 
In \cite{Doganay2019}, two algorithmic approaches are described: the weighted sum method as an example for a scalarization method, see, e.g., \cite{Ehrgott2005}, and a multiobjective descent algorithm as an example for a non-scalarization method, see \cite{Fliege2000}. Gradient descent strategies were implemented for both methods to search for Pareto critical points, i.e., points for which no common descent direction for all objectives exists. 
In this section, we focus on weighted sum scalarizations and present new numerical results for a biobjective test case.

%\paragraph{Weighted Sum Method}
The weighted sum method replaces the multiobjective function $J$ by the weighted sum of the objectives 
$J^{\omega}(\Omega) = \sum_{i=1}^{p}\omega_i J_i(\Omega)$. 
Here, $\omega\geqslant0$ is a weighting vector that represents the relative importance of the individual objective functions. We assume without loss of generality that 
$\sum_{i=1}^p \omega_i=1$. 
The resulting scalar-valued objective function $J^{\omega}$ can then be optimized by (single-objective) gradient descent algorithms, see e.g.\ \cite{CarlGeiger1999}. 
If a global minimum of the weighted sum scalarization $J^{\omega}$ is obtained, then this solution is a Pareto optimal solution of the corresponding multiobjective optimization problem, see, e.g., \cite{Ehrgott2005}. 
The converse is not true in general, i.e., not every Pareto optimal solution can be obtained by the weighted sum method. Indeed, the weighted sum method can not be used to explore non-convex parts of the Pareto front.  
Nevertheless, an approximation of the Pareto front can be obtained by appropriately varying the weights. 

\subsubsection{Case Study and Numerical Implementation}\label{subsec:MOcasestudy}

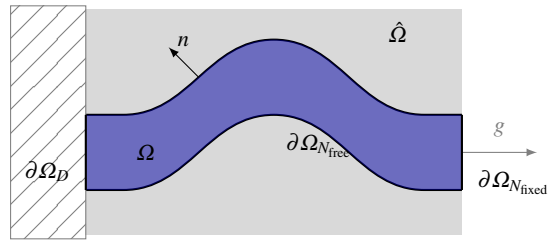
\begin{figure}[ht]
\centering
\begin{tikzpicture}[font={\footnotesize}, scale =1]
			\coordinate (bb_ll) at (0,0.4);
			\coordinate (bb_lr) at (5,0.4);
			\coordinate (bb_tl) at (0,3.4);
			\coordinate (bb_tr) at (5,3.4);
			\coordinate (bar_ll) at (0,1.);
			\coordinate (bar_lr) at (5,1.);
			\coordinate (bar_tl) at (0,2.);
			\coordinate (bar_tr) at (5,2.);
			\fill [fill=black!15!white] (bb_ll) rectangle (bb_tr);
			\draw [thick,name path=A] (bar_ll) to[out=0,in=180] (0.5,1) to[out=0,in=180] (2.5,2.) to[out=0,in=180] (4.5,1) to[out=0,in=180] (bar_lr);
			\draw [thick,name path=B] (bar_tl) to[out=0,in=180] (0.5,2.) to[out=0,in=180] (2.5,3) to[out=0,in=180] (4.5,2.) to[out=0,in=180] (bar_tr);
			\tikzfillbetween[of=A and B] {black!35!blue, opacity=0.5};
			\draw [-latex] (1.5,2.5) -- (1.1,2.9) node [midway, right = 0.1cm, above =0.1cm]{\(n\)} ;	
			\draw [-latex, gray] (5,1.5) -- (6,1.5) node [midway, right = 0.1cm, above =0.1cm]{\(g\)} ;
			\draw[gray, pattern=my north east lines, pattern color= gray, line space=9pt] (0,0.35) rectangle (-1,3.45);
			\draw [thick] (bar_ll) -- (bar_tl) ;
			\draw [thick] (bar_lr) -- (bar_tr) ;	\node[label=left:\(\hat{\Omega}\)] at (4.5,3.1) {};
			\node[label=left:\(\Omega\)] at (1.15,1.5) {};
			\node[label=right:\(\partial\Omega_{N_{\text{fixed}}}\)] at (5,1.1) {};
			\node[label=right:\(\partial\Omega_{N_{\text{free}}}\)] at (2.45,1.6) {};
			\node[label=left:\(\partial\Omega_{D}\)] at (0,1.25) {};
			%\draw[pattern=north east lines, pattern color=black] (0,0.4) rectangle (-0.5,3.4);
			\end{tikzpicture} 
%\sidecaption
\caption{Case study: general setup and starting solution}
\label{fig:case_study_setting}
\end{figure}

In our case study we focus on objectives $J_1$ and $J_2$ as introduced above (i.e., mechanical integrity and cost), and consider a ceramic component $\Omega \subset \R^2$ made from beryllium oxide (BeO)
(with material parameter setting equal to \cite{Doganay2019}, in particular Weibull’s modulus $m = 5$).
%(with Young’s modulus $E = 320\; \text{GPa}$, Poisson’s ratio $\nu = 0.25$, ultimate tensile strength $140\; \text{MPa}$ and Weibull’s modulus $m = 5$)
% calculation of sigma_0 via
% \begin{align}
% \sigma_0 = (\dfrac{\text{uts}^m}{1000})^{\dfrac{1}{m}}
% \end{align}
%
% $m, poisson, E = (2. 0.25\cdot 10^7) $) 
%under tensile load $g= 10^8\,\text{Pa}$ in two dimensions. 
The volume force $f$ is set to $1000\,\text{Pa}$. %\todo{Specify $g$}
See Fig.~\ref{fig:case_study_setting} for an illustration of a possible (non-optimal) shape. 
This shape is used as starting solution for the numerical tests described below. 
The component is of length $0.6\,\text{m}$ and is assumed to have a thickness of $0.1\,\text{m}$. 
It is fixed on the left boundary $\Omega_D$ and the tensile load is acting on the right boundary $\Omega_{N_{\text{fixed}}}$. The parts $\Omega_D$ and $\Omega_{N_{\text{fixed}}}$ are fixed, while the part $\Omega_{N_{\text{free}}}$ can be modified during the optimization process.
The biobjective shape optimization problem is then given by
\begin{equation}
\begin{split}
  \min_{\Omega\in\Oad} & ~J(\Omega):=(J_1(\Omega, u),\ J_2(\Omega))\\
  \text{s.t. } & u \in H^1(\Omega, \R^2) \text{ is the solution of a linear elasticity equation.} %(\ref{stateequation}).
\end{split}\label{ceramicMOP}
\end{equation}
%with $J_1$ concerning the probability of failure, as defined in \eqref{eqa:ObFun}, and $J_2=\int_{\Omega}dx$ concerning the cost.
The component is discretized using a regular $45\times 25 $ grid, see Section~\ref{ShapeOpStruc} and Fig.~\ref{fig:ceramic-gradients}.
% \todo{MB: Bitte gleiche Skalierung für x- und y-Achse und Vektorgrafik verwenden. JS: erledigt}
%\todo{Wenn Platz knapp wird, kann Fig.~\ref{fig:grid+comp} raus.}
%\begin{figure}[h]
%    \centering
%		\includegraphics[width=0.8\textwidth]{mop_figures/grid45x25_movable_boundaries.pdf}
%    \caption{Discretized component; fixed boundary in red, movable boundary in green}
%    \label{fig:grid+comp}
%\end{figure}

\begin{figure}[ht]
	\centering
		\subfloat{\includegraphics[width=0.3\linewidth]{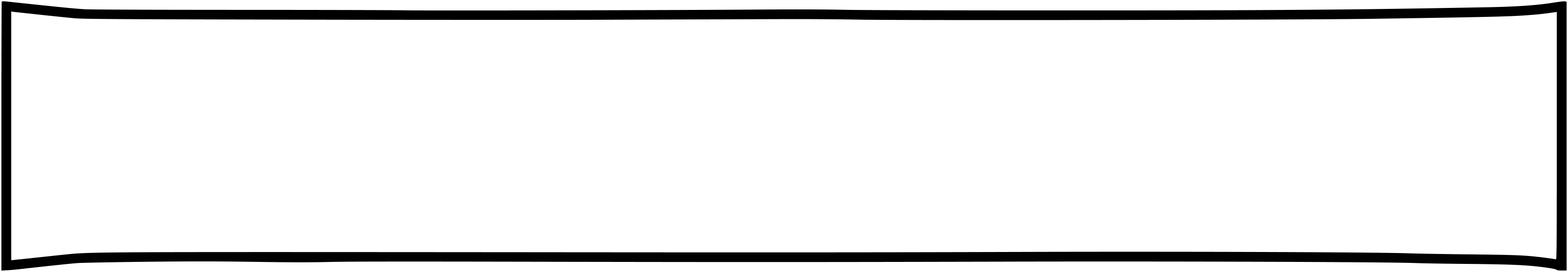}}
        \hfill
		\subfloat{\includegraphics[width=0.3\linewidth]{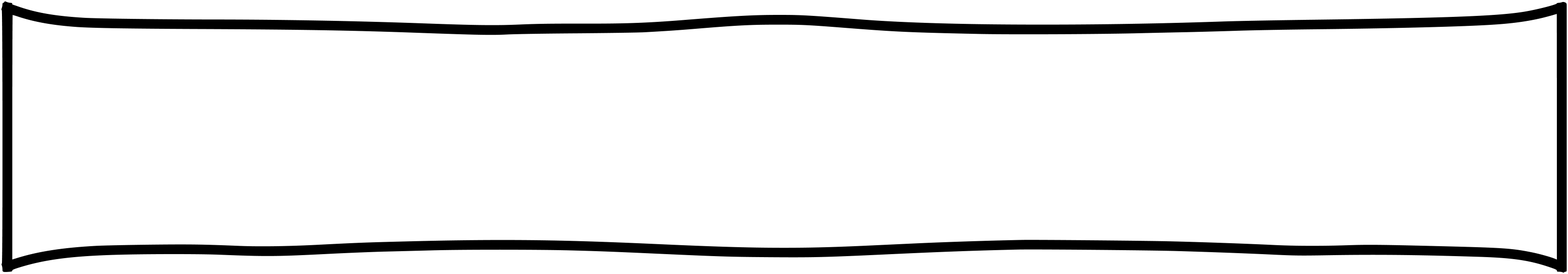}}
		\hfill
		\subfloat{\includegraphics[width=0.3\linewidth]{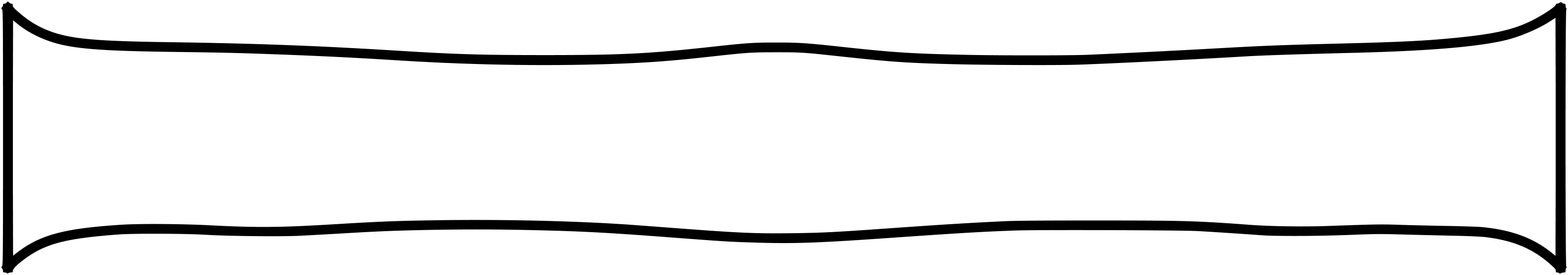}}
	\caption{Near Pareto critical solutions obtained by the weighted sum method}
	\label{fig:wsum}
\end{figure}

\begin{figure}[ht]
% \includegraphics[scale=0.5]{final_iteration133.png}
%\centering
% \begin{tikzpicture} 
% \begin{axis}[axis lines=middle,
% width=7cm,
% xmin=50, 
% xmax=150, 
% ymin=0.03,
% ymax=0.061,
% label style={font=\footnotesize}, 
% tick label style={font=\footnotesize},
% xlabel=$J_1$, % mechanical integrity 
% ylabel={ $J_2$}, %volume
% legend style={cells={anchor=west},font=\scriptsize, at={(axis cs:0.004,120)}}]
% % results of the weighted sum method (with labels)

% \coordinate[style={fill, circle, inner sep = 0.5pt}] (node05) at (axis cs:127.5737072809748298, 0.03739825102054964656) {};
% \coordinate[style={fill, circle, inner sep = 0.5pt}] (node06) at (axis cs:106.9638955229890485, 0.03892833443635850682) {};
% \coordinate[style={fill, circle, inner sep = 0.5pt}] (node07) at (axis cs:89.89224330081458447, 0.04056345323757281490) {};

% \end{axis}
% \end{tikzpicture}

% \centering
\begin{tikzpicture} 
\begin{axis}[axis lines=middle,
width=10cm,
xmin=10, 
xmax=115, 
ymin=0.04,
ymax=0.061,
label style={font=\footnotesize}, 
tick label style={font=\footnotesize},
xlabel=$J_1$, % mechanical integrity 
ylabel=$J_2$, %volume
legend style={cells={anchor=west},font=\scriptsize, at={(axis cs:0.004,120)}}]

% results of the weighted sum method (with labels)
%\node at (axis cs:8.349758844949117531e+01,    4.190074079578058885e-02) [shape=cross out, draw] {}; % w 0.006
%\node at (axis cs:7.476150761149531832e+01,    4.244071576961090386e-02) [shape=cross out, draw] {}; % w 0.007
%\node at (axis cs:7.193553988589690107e+01,    4.280219595026825419e-02) [shape=cross out, draw] {}; % w 0.008

\coordinate (w0006) at (axis cs:8.349758844949117531e+01,    4.190074079578058885e-02); % w 0.006
\coordinate (w0007) at (axis cs:7.476150761149531832e+01,    4.244071576961090386e-02); % w 0.007
\coordinate (w0008) at (axis cs:7.193553988589690107e+01,    4.280219595026825419e-02); % w 0.008
\coordinate (w0009) at (axis cs:6.579354054023639264e+01,    4.369507613091547921e-02); % w 0.009

\coordinate (w0016) at (axis cs:5.198674147491482955e+01,    4.606380141672362233e-02); % w0.016
\coordinate (w0017) at (axis cs:5.092226130310416465e+01,    4.629943197384678299e-02); % w0.017
\coordinate (w0018) at (axis cs:5.109225051871871415e+01,    4.622752877050079201e-02); % w0.018
\coordinate (w0019) at (axis cs:4.697450600945805377e+01,    4.717286336622839027e-02); % w0.019

\coordinate (w0056) at (axis cs:2.272630747873302681e+01,    5.619252956045474312e-02); % w0.056
\coordinate (w0057) at (axis cs:2.243024509908257613e+01,    5.636829403712644015e-02); % w0.057
\coordinate (w0058) at (axis cs:2.223271141326851819e+01,    5.649056747500610237e-02); % w0.058
\coordinate (w0059) at (axis cs:2.205764595811034212e+01,    5.660035980822535756e-02); %w0.059

\coordinate (w002) at (axis cs:4.252814097104339197e+01,    4.833081276286003980e-02); % w0.02
\coordinate (w0025) at (axis cs:3.520734656920360095e+01,    5.059205244100450222e-02); % w0.025
\coordinate (w003) at (axis cs:3.132343335660362271e+01,    5.203204464183875527e-02); % w0.03
\coordinate (w0035) at (axis cs:2.630437211201788728e+01,    5.439839811857595520e-02); % w0.035
\coordinate (w004) at (axis cs:2.726127354757186794e+01,    5.376642965344282515e-02); %w0.04
\coordinate (w0045) at (axis cs:2.447798279811268785e+01,    5.523306610343636791e-02); % w0.045
\coordinate (w005) at (axis cs:2.399630750457910366e+01,    5.545068474828620392e-02); % w0.05
\coordinate (w0055) at (axis cs:2.303005341731193667e+01,    5.600285770318052386e-02); % w0.055
\coordinate (w006) at (axis cs:2.192011206742591156e+01,    5.668770019408336108e-02); % w0.06

%\draw[mark=*, mark options={fill=white}] (w0006) circle (2pt);
%\fill[] (w0007) circle (1pt);
%\fill[] (w0008) circle (1pt);
\draw[mark=*, mark options={fill=white}] (w0006) circle (1pt);
\draw[mark=*, mark options={fill=white}] (w0007) circle (1pt);
\draw[mark=*, mark options={fill=white}] (w0008) circle (1pt);
\draw[mark=*, mark options={fill=white}] (w0009) circle (1pt);

\draw[mark=*, mark options={fill=white}] (w0016) circle (1pt);
\draw[mark=*, mark options={fill=white}] (w0017) circle (1pt);
\draw[mark=*, mark options={fill=white}] (w0018) circle (1pt);
\draw[mark=*, mark options={fill=white}] (w0019) circle (1pt);

\draw[mark=*, mark options={fill=white}] (w0056) circle (1pt);
\draw[mark=*, mark options={fill=white}] (w0057) circle (1pt);
\draw[mark=*, mark options={fill=white}] (w0058) circle (1pt);
\draw[mark=*, mark options={fill=white}] (w0059) circle (1pt);

\draw[mark=*, mark options={fill=white}] (w002) circle (1pt);
\draw[mark=*, mark options={fill=white}] (w0025) circle (1pt);
\draw[mark=*, mark options={fill=white}] (w003) circle (1pt);
\draw[mark=*, mark options={fill=white}] (w0035) circle (1pt);
\draw[mark=*, mark options={fill=white}] (w004) circle (1pt);
\draw[mark=*, mark options={fill=white}] (w0045) circle (1pt);
\draw[mark=*, mark options={fill=white}] (w005) circle (1pt);
\draw[mark=*, mark options={fill=white}] (w0055) circle (1pt);
\draw[mark=*, mark options={fill=white}] (w006) circle (1pt);

\node (contour-w0006) at ($1.3*(w0006)$){\includegraphics[width=1.5cm]{grid36x18_w0_006_final_iteration135.eps}};
%\node (contour-w0007) at ($1.4*(w0007)$){\includegraphics[width=1.5cm]{grid36x18_w0_007_final_iteration151.eps}};
% \node (contour-w0008) at ($1.4*(w0008)$){\includegraphics[width=1.5cm]{grid36x18_w0_008_final_iteration271.eps}};
\node (contour-w0009) at ($1.4*(w0009)$){\includegraphics[width=1.5cm]{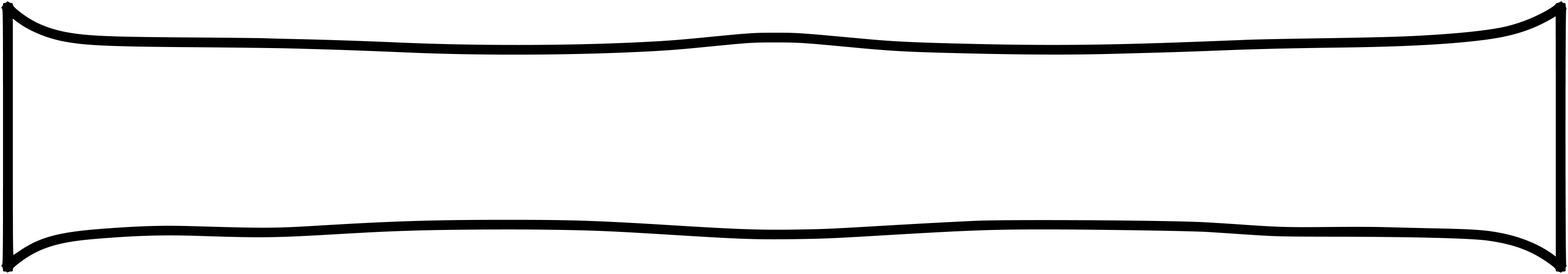}};

% \node (contour-w0016) at ($1.4*(w0016)$){\includegraphics[width=1.5cm]{figs/pareto_front/grid36x18_w0_016_final_iteration233.eps}};
%\node (contour-w0017) at ($1.3*(w0017)$){\includegraphics[width=1.5cm]{figs/pareto_front/grid36x18_w0_017_final_iteration62.eps}};
%\node (contour-w0018) at ($1.4*(w0018)$){\includegraphics[width=1.5cm]{figs/pareto_front/grid36x18_w0_018_final_iteration172.eps}};
\node (contour-w0019) at ($1.3*(w0019)$){\includegraphics[width=1.5cm]{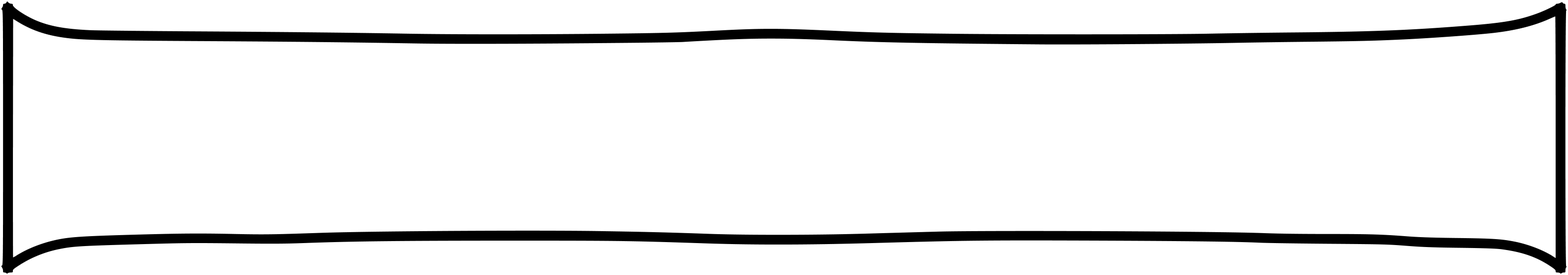}};

%\node (contour-w0056) at ($1.1*(w0056)$){\includegraphics[width=1.5cm]{figs/pareto_front/grid36x18_w0_056_final_iteration106.eps}};
\node (contour-w0057) at ($1.1*(w0057)$){\includegraphics[width=1.5cm]{grid36x18_w0_057_final_iteration465.eps}};

\node (contour-w003) at ($1.3*(w003)$){\includegraphics[width=1.5cm]{grid36x18_w0_030_final_iteration250.eps}};

\draw[very thin, shorten <=2pt] (w0006)-- (contour-w0006) ;
%\draw[very thin, shorten <=2pt] (w0007)-- (contour-w0007) ;
% \draw[very thin, shorten <=2pt] (w0008)-- (contour-w0008) ;
\draw[very thin, shorten <=2pt] (w0009)-- (contour-w0009) ;

% \draw[very thin, shorten <=2pt] (w0016)-- (contour-w0016) ;
%\draw[very thin, shorten <=2pt] (w0017)-- (contour-w0017) ;
%\draw[very thin, shorten <=2pt] (w0018)-- (contour-w0018) ;
\draw[very thin, shorten <=2pt] (w0019)-- (contour-w0019) ;

%\draw[very thin, shorten <=2pt] (w0056)-- (contour-w0056) ;
\draw[very thin, shorten <=2pt] (w0057)-- (contour-w0057) ;

\draw[very thin, shorten <=2pt] (w003)-- (contour-w003) ;

\end{axis}
\end{tikzpicture}
%\sidecaption
     \caption{Approximated Pareto front obtained by the weighted sum method}
    \label{fig:PFapprox}
\end{figure}

We use a gradient descent method to minimize the weighted sum objective function $J^{\omega}$ for different weight vectors $\omega\geqslant0$. 
This is implemented using the negative gradient as search direction and the Armijo-rule to determine a step-size, see e.g.\ \cite{CarlGeiger1999}. 
During the iterations, the component is modified by free form deformations using the method developed in Section~\ref{ShapeOpStruc}. 
Since we have a regular mesh inside the component, only the grid points close to the boundary have to be adapted. 
When the modification during one iteration is too large, a complete remeshing is performed, still using the approach described in Section~\ref{ShapeOpStruc}. 
%\todo{MB: Trat das wirklich auf?}
To avoid oscillating boundaries and overfitting, we apply a regularization approach based on \cite{Schulz2016}. 
%
%\paragraph{Results} %Test Case}
Numerical results for three choices of the weight vector $\omega$ are shown in Fig.~\ref{fig:wsum}, %\todo{Welche Gewichte? Entsprechend sortieren. JS: gebe hier keine gewichte an, da die Zielfunktionen vorskaliert sind}
%\todo{Bei knappem Platz kann Fig.~\ref{fig:wsum} raus.}
and an approximation of the Pareto front is given in Fig.~\ref{fig:PFapprox}. In these cases no remeshing step had to be performed, because the step length was restricted to the mesh size. 
%\todo{Skalierung der Achsen? JS: Bild ist noch in Arbeit}
%\todo{Fig.~\ref{fig:grid+comp} und \ref{fig:wsum} können bei Platzmangel raus}

%%%%%%%%%%%%%%%%%%%%%%%%%%%%
\subsubsection{Gradient Enhanced Kriging for Efficient Objective Function Approximation}
%\paragraph{Gradient Enhanced Kriging for Efficient Objective Function Approximation} %\todo{JS: Eigene Section?}
To cut computational time of the optimization process one can apply surrogate models to estimate expensive to compute objective functions.
Optimization on the surrogate model is relatively cheap and yields new points which then in a next step are evaluated with the expensive original objective function. In the biobjective model presented above, the mechanical integrity $J_1$ is expensive, while the volume $J_2$ can be easily evaluated. 
We thus suggest to replace only the expensive objective $J_1$ by a model function.

Let $\{\Omega_1,\hdots, \Omega_M\} \subset \Oad$ be sampled shapes with responses $\{y_1,\hdots, y_M\}:=\{J_1(\Omega_1),\hdots, J_1(\Omega_M)\}$. 
\emph{Kriging} is a type of surrogate model that assumes that the responses $\{y_1,\hdots, y_M\}$ are realizations of Gaussian random variables $\{Y_1,\hdots, Y_M\}:=\{ Y(\Omega_1),\hdots,  Y(\Omega_M)\}$ from a Gaussian random field $\{Y(\hat{\Omega})\}_{\hat{\Omega} \in \Oad}$. For an unknown shape $\Omega_0$ the Kriging model then predicts 
\begin{equation*}
    \hat{y}(\Omega_0)=\mathbb{E}[Y(\Omega_0)| Y(\Omega_1)=y_1,\hdots, Y(\Omega_M)=y_M],
\end{equation*}
i.e., the estimated objective value of $\Omega_0$ is the conditional expectation of $Y(\Omega_0)$ under the condition that the random field is equal to the responses at the sampled shapes, or in other words the predictor is an interpolator. 
An advantage of this method is, that the model also provides information about the uncertainty of the prediction, denoted as $\hat{s}(\Omega_0)$, see \cite{Keane2008} for more details.

If, as in our case, gradient information is available, one can incorporate this into the Kriging model which is then called \emph{gradient enhanced Kriging}. 
One follows the same idea: the gradients $\{\Bar{y}_1,\hdots,\Bar{y}_{\Bar{M}}\}:=\{ \nabla J_1(\Omega_1),\hdots, \nabla J_1(\Omega_{\Bar{M}})\}$ are assumed to be realizations of the Gaussian random variables/vectors $\{\Bar{Y}_1,\hdots,\Bar{Y}_{\Bar{M}}\}$. 
Adding these random variables to the ones w.r.t.\ the objective values enables one to predict objective values and gradients at unknown shapes 
$\Omega_0$, see also \cite{Keane2008} for more details.

In the optimization choosing the predictor $ \hat{y}(\Omega_0)$ as the objective to acquire new points to evaluate with the original function may yield poor results. 
Since then one assumes that the prediction has no uncertainty, i.e.\ $\hat{s}(\Omega_0)=0$, and areas, for which the predictor has bad values and a high uncertainty while the original function has better values than the best value at the moment, may be overlooked. 
Hence, one has to choose an acquisition function that incorporates $\hat{y}(\Omega_0)$ and $\hat{s}(\Omega_0)$ to balance the exploitation and exploration in the optimization.

We note that this gradient enhanced Kriging approach is a direction of ongoing development for the in house optimization process {\tt AuoOpti} at the German Aerospace Center (DLR), see e.g.\ \cite{schaefer2012multiobjective,kroger2010axisymmetric} for design studies using the {\tt AutoOpti} framework. In our future work, we therefore intend to benchmark the EGO-based use of gradient information with the multi criteria descent algorithms and identify their respective advantages for gas turbine design.

%\todo{MB: Kriging wird derzeit (noch) nicht verwendet, oder? Warum brauchen wir dann diesen Abschnitt?}

%%%%%%%%%%%%%%%%%%%%%%%%%%%%%%%%%%%%%%%%%%  (Lucas Maede & Benedikt Engel)

%%%%%%%%%%%%%%%%%%%%%%%%%%%%%%%%%%%%%%%%%%%%%%%%%%%%%%
\section{Applications}\label{sec:applications}

In this section we present the industrial implications of the GivEn consortium's research. The German Aerospace Center (DLR), Institute for Propulsion Technology, here is an important partner with an own tool development that involves an in-house adjoint computational fluid dynamics solver 
%\todo{Auf TRACE und AutoOpti wird später im Kapitel nicht mehr eingegangen - sind wirklich keine Rückkopplungen von GIVEN-Ergebnissen auf die beiden Softwarepakete vorgesehen?}
% Anmerkung: [J.B.] Die Rückwirkung wird eher über Prozesse um TRACE und adjoint TRACE herum entstehen. AutoOpti hat eine sehr starke Fokussierung auf Ersatzmodellbasierte Optimierung weshalb eine Implementierung dort nicht einfach sein wird. Es könnte eher auf eine Kopplung von Gradientenabstieg mit AutoOpti hinauslaufen, indem Lösungsvorschläge aus Gradientenabstieg als neu vorgeschlagene Lösung mit in die nächsten Ersatzmodelle einfließen.
TRACE as well as a multi-criteria optimization toolbox AutoOpti. 
As TRACE and AutoOpti are widely used in the German turbo-machinery industry, a spill-over of GivEn's method to the DLR assures an optimal and sustainable distribution of the research results.

In a second contribution, Siemens Energy shows that results developed in the GivEn project can also be directly used in an industrial context, taking multi-criteria tolerance design as an example.  

%%%%%%%%%%%%%%%%%%%%%%%%%%%%%%%%%%%%%%%%%% DLR
\subsection{German Aerospace Center (DLR)}\label{sec:31} 
Industrial turbomachinery research at DLR includes several activities that benefits from the insights gained in this project.
These activities primarily pursue two goals:
\begin{enumerate}[(i)]
   \item\label{dlrAppItem1} Assessment of technology potential for future innovations in the gas turbine industry.
   \item Development of efficient design and optimization tools that can be used, for instance, to perform \eqref{dlrAppItem1}.
\end{enumerate}

\subsubsection{Challenges of Industrial Turbomachinery Optimization}
The gain in aerodynamic performance of both stationary gas turbines and aircraft engines that has been achieved over the last decades, leaves little room for improvement if solely aerodynamics is considered. 
More precisely, aerodynamic performance enhancements that neglect the issues of manufacturing, structural dynamics or thermal loads, will typically not find their way into application. One of the reasons for this is the fact that real engines currently designed already have small "safety" margins. 
Summarizing, one can conclude that aerodynamic performance, structural integrity as well as manufacturing and maintenance costs have become competing design goals. 
Therefore, the design of industrial turbomachinery has become a multi-disciplinary multi-criteria optimization problem. Accordingly, DLR is highly interested in advances concerning both simulation tools for coupled problems and multi-disciplinary optimization (MDO) techniques. 

Multi-criteria optimizations based on high-fidelity simulations are firmly established in the design of turbomachinery components in research and industry \cite{thevenin2008, li2017}.
As explained above, current developments increasingly demand the tighter coupling of simulations from multiple disciplines. Reliable evaluations of such effects require the simultaneous consideration of aerodynamics, aeroelasticity and aerothermodynamics. 
Moreover, optimization should account for the influence of results from these disciplines on component life-times.

Gradient-free optimization methods, typically assisted by surrogate modeling, prevail in today's practical design processes \cite{peter2008, samad2009}.
A tendency towards optimizing with higher level of detail and optimizing multiple stages simultaneously leads to higher dimensional design spaces, making gradient-free methods increasingly expensive and gradient-based optimization the better suited approach.

\subsubsection{Expected Impact of GivEn Results}
The methods described in the preceding sections describe how derivatives for a fully coupled aerothermal design evaluation process can be computed efficiently and how a gradient based optimization procedure can be constituted for the design criteria of efficiency and component life-time.
The exemplary process, developed in the frame of this project, serves as a research tool and a base to adopt the methods for other applications resulting in different levels of simulation fidelity, different sets of disciplines as well as different objectives.
The partners from DLR consider the goals of this project as an important milestone that could enable researchers to tackle, among others, the following problems:
\begin{itemize}
    \item Concurrent optimization of turbine aerodynamics and cycle performance with the goal of reducing cooling air mass flows. Such optimization should take into account the redistribution and mixing of hot and cold streaks to be able to predict aerodynamic loads in downstream stages.
    \item Assessment of the technology potential of operating the burner at partial admission (or even partial shutdown) conditions in order to achieve good partial load performance while avoiding a significant increase in emissions. 
    \item Assessment of potentials of thermal clocking. The idea is to be able to reduce cooling air if the relatively cold wake behind cooled vanes is used in downstream stages. Such optimizations will be based on unsteady flow predictions.
\end{itemize}
The necessary changes to create adjoints for existing evaluation processes stretches from parametrization to the simulation and post-processing codes for the different disciplines involved. Moreover, DLR expects to benefit from 
the coupling strategies developed here. These should be both sufficiently accurate and apt for an appropriate strategy to define coupled adjoint solvers. A particularly important milestone that this project is to achieve is the establishment of a life-time prediction that is suitable for gradient-based optimization. 

DLR will not only apply these advances in aerothermal design problem but also hopes that general conclusions can be drawn that carry over to other multi-disciplinary design problems that involve coupled simulations of flows and structures.

Theoretical concepts are shared in joint seminars in the early stages of the project to lay the foundation of common implementation of prototypes in later phases of the project.
The role of these prototypes is to explore the implementation of the methods and simultaneously are used to communicate about changes to existing design evaluation chains and requirements for practical optimizations.

%%%%%%%%%%%%%%%%%%%%%%%%%%%%%%%%%%%%%%%%%% Siemens
\subsection{Siemens Energy}\label{sec:32} 

In accordance with the final aim of the GivEn project, the main usage of adjoint codes in the turbomachinery industry concerns the development of optimization process tool chain for design purpose.
%To reach this goal, a substantial amount of work must be invested in research and development delaying the beneficial impact of adjoint codes for the industry. (removed by request of EARS Reviewer)
Recently the use of adjoint codes to consider the effects of manufacturing variations has been proposed.  This new application can have a major impact in the very competitive market of gas turbine for power generation at relatively low cost and time horizon. The following sections describe some aspects of the development of such a method using the example of a turbine vane.

%%%%%%%%%%%%%%%%%%%%%%%%%%%%%%%%%%%%%%%%%%%%%%%%%%%%%%%%%%%%
\subsubsection{Effects of Manufacturing Variations}
The vanes and blades of a gas turbine are designed to deflect and expand the flow leaving the combustion chamber and generate torque at the rotor shaft.
It implies that these components must operate at very high temperature while maintaining good aerodynamics and mechanical properties. 
Therefore, the material chosen for such hot gas section parts usually belongs to the Ni-base superalloy family and the manufacturing process involves precision casting which is a very complex and expensive technique.
% Rewritten according to EARS reviewer request
During the manufacturing process, the shape accuracy of each casting is assessed by geometrical measurements.
%During the manufacturing process, the quality of a vane is assessed by geometrical measurements made at several well-chosen locations. 
If all of the measured coordinates lie within a specified tolerance band, the part is declared compliant and otherwise scrapped. 
%This quality assessment suffers from several drawbacks: it is relatively costly and time consuming and the number of points investigated - based on experience - is very limited . %This rather crude process leads to a relatively high scrap rate which increase the price of a turbine.
In order to reduce uneconomic scrapping, the re-evaluation of acceptable tolerances is a constantly present challenge. 
In the past, point-based measurements and subsequent combination with expert judgement were conducted.
Nowadays, advanced 3D scanners enable the acquisition of highly detailed geometric views of each part in a short amount of time which are better suited to characterize the analyze geometrical deviations.
%Despite the fact that the control points are wisely chosen it would be possible that some geometrical deviations elsewhere on the vane might have a positive or negative influence on the vane's characteristics.
They furthermore allow to focus on the $\it{effects}$ of geometric variations on the component's functions rather than on the geometric variations themselves.
%It would be therefore more interesting to use a more detailed geometrical description of the manufactured vanes and to focus on the $\it{effects}$  of the geometric variations on the vane's functions rather than on the geometric variations themselves.
%The recent development of 3D scanner enables to obtain an highly detailed geometric view of each part in a short amount of time.
Nonetheless, evaluating the characteristics of each component produced using traditional computational methods - CFD and FEM - would not be industrially feasible due to the computational resource required.
In this context, the capability of adjoint codes to take into account very efficiently the effects of a great number of parameters on one objective function is a crucial asset.
It allows to benefit from the high resolution scans and model the manufacturing variations on the surface directly as the displacement of each point on the surface. 

However, since adjoint equations are linearized with respect to an objective function, the accuracy of the predictions provided with the help of an adjoint code will be limited by the magnitude of the deformations.
In other words, the deformations must be small enough such that a linear approximation of the effects on the objective function is appropriate.
Therefore the usefulness of adjoint codes in an industrial context must be investigated using $\it{real}$ manufacturing deviations. To this aim, 102 scans of heavy-duty turbine vanes have been used for the validation of the adjoint codes.
Since the role of a turbine vane consists in fulfilling aerodynamic and mechanical functions, the objectives chosen for the investigation must be chosen within these two families.

%%%%%%%%%%%%%%%%%%%%%%%%%%%%%%%%%%%%%%%%%%%
\subsubsection{Validation of the Tools}
In a context of a more volatile energy production induced by renewable sources, gas turbines must be started or stopped very often and produce energy the most efficiently possible. 

\begin{figure}[t]
	\centering
		\subfloat{\includegraphics[width=0.45\linewidth]{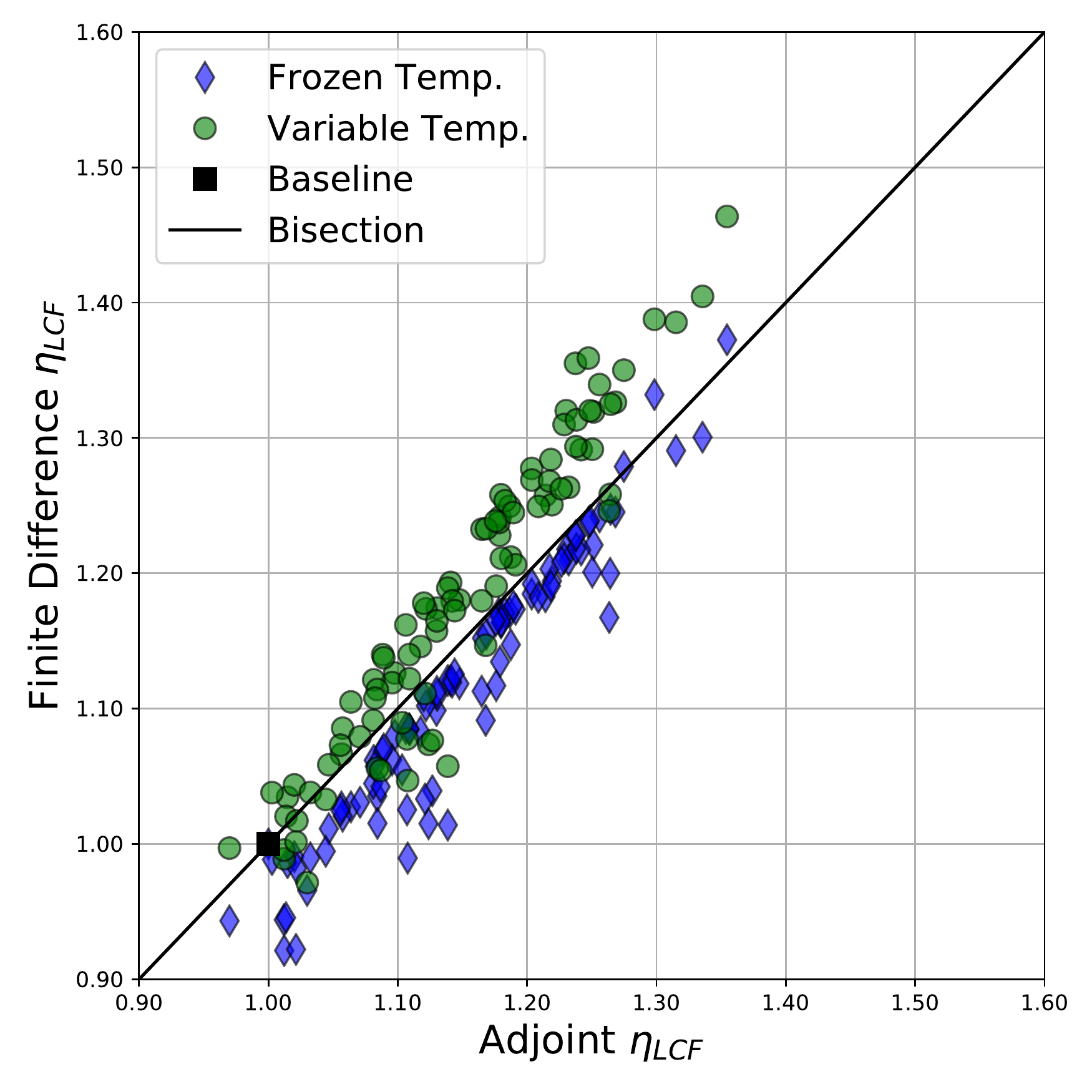}}
		\hfill
		\subfloat{\includegraphics[width=0.45\linewidth]{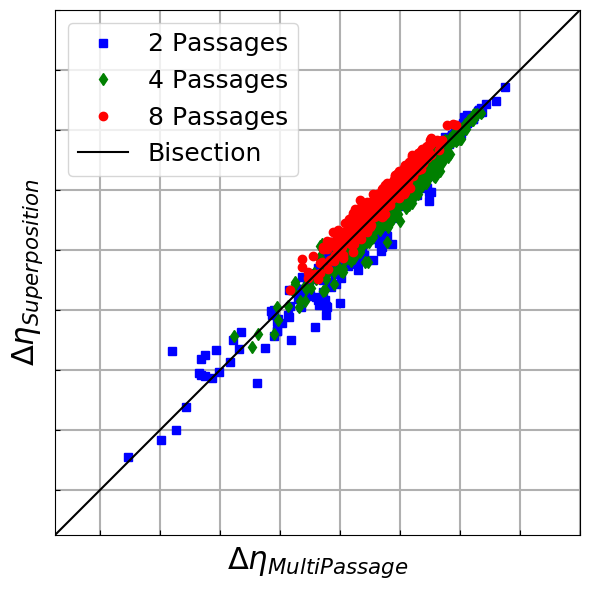}}
	\caption{Weibull scale parameter ($\eta_{LCF}$, right) and isentropic efficiency ($\eta$, left) obtained by adjoint codes versus finite differences.}
	\label{fig:SIEMENS_GRADIENT}
\end{figure}

% LCF
For turbine vanes, starting or stopping a machine more often implies that important temperature gradients will be experienced also more often leading potentially to LCF problems.
The occurrence of LCF issues can be numerically estimated by computing the maximal number of start/stop cycles above which a crack appears on the vane surface. For this assessment, the probabilistic Weibull-based approach by Schmitz \textit{et al.} has been chosen.
It was combined with the primal and adjoint equations and has been implemented into a Siemens in-house software.
In Fig.~\ref{fig:SIEMENS_GRADIENT}, left, a comparison of the gradients of the Weibull scale parameter $\eta_{LCF}$ computed with the adjoint code and the finite difference for each of the 102 vanes considered is presented. Actually two sets of data are present: one assuming that the temperature field on the vane is independent of the actual geometry - the so-called frozen temperature hypothesis - while the other considers on the contrary a variable temperature field. In both cases, all points are very close to the first bisectional curve, meaning that for each of the 102 vanes, the manufacturing deviations have a linear impact on the Weibull scale parameter. In other words, the magnitude of the real manufacturing deviations are small enough so their effects on LCF can be successfully taken into account by an adjoint code.   
A more detailed investigation can be found in the publication of Liefke \textit{et al.}~\cite{liefke2020}.

% Aero
The efficiency of a gas turbine depends of course on the aerodynamic characteristic of the vanes. However the behaviour of each single vane makes only sense if they are integrated into a row together. 
When considering manufacturing variations, an additional difficulty appears namely the interaction of the manufacturing deviations of adjacent blades with each other. In this case, it would be interesting to compare the gradients obtained with the adjoint code versus the finite difference, not only for a purely axis-symmetrical situation (e.g.\ one passage) but also for several different arrangements (e.g.\ more than one passage). 
Fig.~\ref{fig:SIEMENS_GRADIENT}, right, presents the gradients of the stage isentropic efficiency for such a configuration using DLR's CFD suite TRACE ~\cite{backhaus2017application}. It can be seen that independently of the arrangement considered - 2,4 or 8 passages - there is a very good agreement between the gradients predicted by the adjoint code and those obtained with the help of the finite difference. The work of Liefke \textit{et al.}~\cite{liefke2019} summarizes the complete investigation of this case.

%%%%%%%%%%%%%%%%%%%%%%%%%%%%%%%%%%%%%%%%%%%%%%
\subsubsection{Industrial perspectives}
%\todo{Im gesamten Kapitel 3.2 kommt GiVen nur im ersten Satz vor. Bei der Diskussion der drei Teilkapitel wäre es schön, wenn ab und zu ein Hinweis kommt, wie denn GiVen zu den beschriebenen  Aufgaben von Siemens Energy beitragen kann. Ich denke, nur immer von Adjoint approach zu sprechen, ist da zu wenig.}
The previous section demonstrated that adjoint codes can be successfully used in an industrial context to quantify the impact of manufacturing variations on low-cycle fatigue and isentropic efficiency. 
It could possible to extend this approach to additional physical phenomena such 
as high-cycle fatigue or creep, given that the equations modelling these phenomena can be 
differentiated and of course that the impact of the manufacturing variation remains linear.

In addition to the direct and short-term benefits for the manufacturing of turbomachine components, the results presented in this section also demonstrate that adjoint codes can be successfully deployed and used in an industrial context. The confidence and experience gained will pave the way for other new usage within Siemens Energy.
Especially, the tools and concepts developed within the GiVen project will greatly contribute to the creation of more rapid, efficient and robust multidisciplinary design optimization of compressors and turbines either based fully on adjoint codes or in combination with surrogate models.

%%%%%%%%%%%%%%%%%%%%%%%%%%%%%%%% ACKLOWLEDGEMENT %%%%%%%%%%%%
\begin{acknowledgement}
The authors were partially supported by the
BMBF collaborative research project GivEn under the grant no.\ 05M18PXA.
\end{acknowledgement}

%%%%%%%%%%%%%%%%%%%%%%%%%%%%%% references %%%%%%%%%%%%%%

%%%%%%%%%%%%%%%%%%%%%%%% Springer-Verlag %%%%%%%%%%%%%%%%%%%%%%%%%%
%
% BibTeX users please use
% \bibliographystyle{plain}
% \bibliography{literature.bib}
%

%%%%%%%%%%%%%%%%%%%%%%%%%%%%%%% end of GIVEN paper %%%%%%%%%%%
\end{document}